\input amssym.def
\input amssym.tex

\def\today{\ifcase\month\or
     January\or February\or March\or April\or May\or June\or
     July\or August\or September\or October\or November\or December\fi
     \space\number\day, \number\year}

\def\newline{\hfil\break}

\def\u{\underline}
\def\o{\overline}

\def\Spf{{\rm Spf\,}}
\def\id{{\rm id}}
\def\Fil{{\sl Fil}}

\def\leq{\leqslant}
\def\geq{\geqslant}


\def\ds{\displaystyle}

\def\ov{\overline}
\def\r{\right}
\def\l{\left}

\def\phi{\varphi}

\def\theta{\vartheta}


\font \title=cmbx10 scaled \magstep5

\font \tpar=cmbx10 scaled \magstep1
\font \tparit=cmmi10 scaled \magstep1


\def\b#1{ {\Bbb {#1}} }

\def\c#1{ {\cal #1} }





\def\aqq{\left[\!\left[}
\def\cqq{\right]\!\right]}

\newcount\notenumber
\def\clearnotenumber{\notenumber =0}
\def\note{\global\advance\notenumber by 1 \footnote
              {$^{\the\notenumber}$} }

\newbox\corrente

\newcount\n
\def\clearn{\n =0}
\def\N{\global\advance\n by 1
             \global\clearnn
             \global\clearnnn
             \global\clearnnnn
             \global\clearnotenumber
             \global\setbox\corrente=\hbox{{\the\n}}
             {\the\n}.
           }

\newcount\nn
\def\clearnn{\nn =0}
\def\NN{\global\advance\nn by 1
              \global\clearnnn
              \global\clearnnnn
              \global\setbox\corrente=\hbox{{\the\n}.{\the\nn}}
              {\the\n}.{\the\nn}%
           }

\newcount\nnn
\def\clearnnn{\nnn =0}
\def\NNN{\global\advance\nnn by 1
              \global\clearnnnn
              \global\setbox\corrente=\hbox{{\the\n}.{\the\nn}.{\the\nnn}}
              {\the\n}.{\the\nn}.{\the\nnn}%
           }

\newcount\nnnn
\def\clearnnnn{\nnnn =0}
\def\NNNN{\global\advance\nnnn by 1
 
\global\setbox\corrente=\hbox{{\the\n}.{\the\nn}.{\the\nnn}.{\the\nnn}}
              {\the\n}.{\the\nn}.{\the\nnn}.{\the\nnn}%
           }


\clearn
\clearnotenumber

\font \proc=cmcsc10
\def\proclaim#1. #2\par
        {\medbreak  {\proc #1. \enspace }{\sl #2\par }
          \ifdim \lastskip <\medskipamount \removelastskip
                  \penalty 55\medskip \fi}

\def\biblio#1#2{ \item{\hbox to 1cm{[#1]}}
                 {#2} }

\def\epsilon{\varepsilon}
\def\G{\Gamma}
\def\Gp{\Gamma_p}
\def\gp{\gamma_p}
\def\m{(-)}

\def\mod{{\rm mod}}


\def\NoBlackBoxes{\global\overfullrule 0pt}
\NoBlackBoxes

\headline={\hfil}
\footline={\hfil {\rm\folio} \hfil}
\hsize=15 true cm

\centerline{\tpar {\tparit p}-{adic formulas and unit root {\tparit F}-subcrystals}}
\centerline{\tpar {of the hypergeometric system}
\footnote{}{\tenrm 
Version of September 7, 2004}}
\footnote {}{AMS Subject
Classification:11T23,11S31,12H25,14F30.}

\vskip 10 pt
\centerline{\rm by }
\vskip 10 pt
\centerline{\proc Francesco Baldassarri {\rm and} Maurizio Cailotto}
\vskip 30 pt

\n =-1


{\bf \N{Introduction.}}
\vskip 5 pt

This article is dedicated
to the study of $p$-adic analytic continuation of the unit root
$F$-subcrystal of a logarithmic $F$-crystal in
the open tube of a
singularity. We will show how
the possibility of extending that unit root
crystal as a {\it non-singular} crystal in more singular classes,
originates non-trivial formulas
of analytic continuation
of classical functions. As Katz puts it in [Ka], this should be
considered a formula of analytic continuation {\it par excellence}.
We improve on Katz' treatment at least in  that we allow  logarithmic
singularities. A typical example of a formula of this type is the
Koblitz-Diamond formula [Ko], [D], for the analytic
continuation of the function
${\cal F}(a,b,c; \lambda)$ [$p$-DE IV] related to
the classical Gauss hypergeometric function
$F(a,b,c;\lambda)$,
for $a$, $b$, $c \in {\b Z}_p$, $c \notin {\b Z}_{\leq 0}$.
We recall 
that, for any $(a,b,c) \in ({\b  Z}_p)^3$ for which it
makes sense,
$$F(a,b,c;\lambda) = \sum_{s=0}^{+\infty}{{(a)_s(b)_s}\over{(c)_ss!}}
\lambda^s
\in {\b Q}_p[[\lambda]] \; ,$$
where, as usual, we use Pochammer's notation $(a)_s = a(a+1) \dots (a+s-1)$.
The function ${\cal F}(a,b,c; \lambda)$ is the maximal $p$-adic 
analytic extension of the
ratio
$$ {
{F (a,b,c; \lambda)}
\over
{F (a^{\prime},b^{\prime},c^{\prime}; \lambda^p)}
}\in 1+\lambda \b Q_p \aqq \lambda \cqq
$$
where for $a \in {\b Z}_p$, $a^{\prime}\in {\b Z}_p$ is
uniquely
defined by the condition that
$pa^{\prime} - a = \mu_a \in \{0,1,\dots,p-1\}$
(We also recursively define $a^{(0)} = a$,
and $a^{(i+1)} = (a^{(i)})^{\prime}$, for $i = 0,1, \dots$).
The Koblitz-Diamond formula asserts that if $c^{(i)} \in {\b Z}_p^{\times}$ and
$ \mu_{c^{(i)}} \geq \mu_{a^{(i)}} + \mu_{b^{(i)}}$
for any $i = 0,1,\dots$, then
${\cal F}(a,b,c; \lambda)$ extends analytically to
the open disk of radius 1 around $\lambda=1$ and
$${\cal F}(a,b,c;1) = {
{\Gamma_p(c)
\Gamma_p(c-a-b)
}
\over
{\Gamma_p(c-a)
\Gamma_p(c-b)}
}
$$
(where $\Gamma_p$ denotes the
Morita $p$-adic gamma function),
see [D], and [Ko] for $c=1$.
A similar discussion, in the
non-singular case, appears in  [Ba], where it is used  to explain a
formula of Young [Yo]:
If  $\mu_{a^{(i)}} \leq \mu_{b^{(i)}}< p-1$, $\mu_{a^{(i)}}$ even, and
$2 \mu_{b^{(i)}} - \mu_{a^{(i)}} \leq p-1$
for all $i \in \b N$
then $\c F (a,b,1+a-b;\lambda )$ admits an analytic extension
to the class of $-1$, and
$$
     \c F (a,b,1+a-b;-1 ) = \m^{\mu_a \over 2}
    {{\Gp \l( {a \over 2}\r) \Gp \l( b - {a \over 2} \r) }
     \over { \Gp (a) \Gp (b-a) }} \; .
$$
The original paper by Young needs the further condition  that $a,b$
be rational (so in ${1 \over {p^f -1}} \b Z$, for some $f\in\b Z_{\geq 1}$).  The
generalization and the interpretation in terms of unit-root
$F$-subcrystal is given in [Ba].
\par
The reader may have noticed already that our discussion goes
beyond the classical theory [Ka] of $F$-crystals, even if extended to
the framework of logarithmic schemes [Sh]. What we  are really
dealing with is a {\it Dwork family of (filtered, logarithmic)
$F$-crystals}: a structure which was introduced by the first author
in lectures at a $p$-adic Summer School held in Trento in June 1995,
inspired on Dwork's treatment of hypergeometric differential modules
[GHF]. Since those lectures have  unfortunately   remained, as of
today,  unpublished, we quickly define Dwork families of
$F$-crystals in section 2  below, with the promise of making
available a more satisfactory treatment as soon as possible. The
main point  is that  our (logarithmic) crystals $M(a,b,c)$ on, say,  the
formal $p$-adic base $\hat {X}$ = $\hat {\b P}$ := $p$-adic
completion of ${\b P}^1_{\b Z_p}$, with the log-structure induced by
the three  ${\b Z_p}$-points $\{ 0,1,\infty \}$ (or   $ \hat {X} = {\rm Spf} \
{\b Z_p} \{ \lambda , {1 \over {\lambda (1 - \lambda)}} \}$, if one
prefers to avoid singularities) depend on parameters $a$, $b$, $c$ in
$\b Z_p$, and Frobenius is a horizontal transformation  $F(\phi) :
\phi^{\ast}M(a^{\prime},b^{\prime},c^{\prime}) \longrightarrow
M(a,b,c)$, where $\phi$ is some lifting of the absolute Frobenius to
(some open subset of ) $\hat{X}$. This is reminiscent of  Mazur's
theory of $F$-spans [Mz], but differs from it in one crucial point.
The underlying ${\cal O}_{\hat {X}}$-module  of $M(a,b,c)$ is here
independent of $a,b,c$, while the connection $\nabla_{a,b,c}$ varies
with the parameters. In that respect, the hypergeometric  $M(a,b,c)$
is simply $ {\cal O}_{\hat {\b P}}^2$, and we will always express the
connection and the Frobenius map in terms of this trivialization.
This choice of global basis for the hypergeometric module, is also
compatible with the filtration, which however varies with $(a,b,c)$
(morally, with $(a,b,c)$ mod $p$).
   When $a,b,c$ are rational
numbers, say in
$(p^f -1)^{-1} \b Z$, after $f$ steps Frobenius gets
us back to $M(a,b,c)$, and an $f$-th iterate of Frobenius becomes the
standard semilinear automorphism  of classical $F$-isocrystals, with
respect to $\phi^f$.
\par
The notion of a Dwork family of
$F$-crystals is actually taken to be more flexible than what we just
said. According to Dwork's taste, we want to allow for integral
translations of the parameters $a$, $b$, $c$, to reflect the
existence of Gauss' contiguity relations on classical hypergeometric
functions [Po], [GHF], [Bo], [Ku], [Ba]. These relations can be used
in great generality to determine the finite-difference equation
behavior of the Frobenius matrix, and they determine its shape, at
least under the assumption that that matrix be $p$-adic meromorphic
as a function of the parameters $a$, $b$, $c$ in $\b Z_p$. This
assumption is known as the {\it Boyarsky principle}, and is know to
hold for hypergeometric functions  [GHF, 4.7.1]. The matrix
$\gamma^{(\phi)} (\vec a, \vec b ; \lambda)$  used in  this article
is the one of [Bo, 3.1] and [Ku, 1.12.3]. (When $\vec a = (a,b,c) \in
{\b Z}_p^3$, $\vec b = (a^{\prime},b^{\prime},c^{\prime})$ and $\phi
(\lambda) = \lambda^p$, it coincides with the matrix
$\left({B_1\atop B_3}{B_2\atop B_4}\right)$
of [LDE, 4.5.1]). The precise relation
between the Frobenius matrix described in [GHF, 4.7.1] and the one of
[LDE, 4.5.1] is given in [Ba, 2.21].  We use this finite-difference
method  in the appendix  to provide an alternative calculation of the
dominant polar term of the Frobenius matrix in the singular classes
at 0 and 1.
\par

The main  assumption in this paper is the
splitting of the Frobenius matrix of the hypergeometric family
in  the following cases [LDE, 6.6]:
\itemitem{(1)} $\mu_c < \min \l( \mu_a , \mu_b \r) $,
\itemitem{(2)} $\mu_c > \max \l( \mu_a , \mu_b \r) $.
\par \noindent
If $a,b,c$ are rational, this type of condition leads to the
existence of a
   unit root $F$-subcrystal
(of rank one) of the logarithmic $F$-crystal of rank two associated
to the hypergeometric system.

\par
The $p$-adic theory of the hypergeometric system
$${d \over {d\lambda}}Y = Y \pmatrix{ -{c \over \lambda} &
{{c-a} \over {1-\lambda}}\cr {{c-b} \over \lambda} &
{{a+b-c} \over {1-\lambda}}\cr} $$
has been deeply investigated by Dwork [LDE].
The explanation of the Koblitz-Diamond formula
preliminarly requires the calculation of
the eigenvalues of Frobenius operating on the
eigenvectors of classical monodromy
in the class of a single logarithmic singularity. For this
calculation to make sense, we must ensure the convergence of the
uniform part of the classical  fundamental solution matrix at the
singular points  in an open  disk of radius 1. This forces us into a
rather bizarre domain for $(a,b,c) \in {\b Z}_p^3$, ensuring that all
exponents and exponent differences of our differential system consist
of $p$-adically non-Liouville numbers. The standard transfer theorems
of $p$-adic analysis can then be applied [Ch], [BC2], [DGS, Chap. 6].
For the hypergeometric system,
the singular class of 0, and with further restrictions on $(a,b,c)$
({\it e.g.} that they be in $({\b Z}_p \cap {\b Q})^3$), this is
the difficult computation of
chapters 24 to 26 of [LDE].

To transfer that information to the class of, say, 1,
one must understand sufficiently well
the action of M\"obius transformations on
the solutions of the hypergeometric equation.
This action has also been analysed by Dwork in [Ku], a remarkable
paper that adds also to the classical cohomological
understanding of the Kummer transformations, and
determines the effect of those transformations
on the Frobenius matrix. We need to complement section 4 of [Ku]
with a more flexible formula for the changes in the
Frobenius matrix.
   \par For simplicity, we will  further assume in
our calculations that the monodromy at 0 and 1 is semisimple as in
[LDE, Chapter 25]: this assumption does not affect the final result.
One can in fact similarly extend the results of [LDE, Chapter 26], to
cover the case of logarithmic solutions.  In any case, our method
proves the Koblitz-Diamond formula under the slightly modified
assumptions  $\min (\mu_{a^{(i)}}, \mu_{b^{(i)}}) > 0$ and
$ \mu_{c^{(i)}} > \mu_{a^{(i)}} + \mu_{b^{(i)}}$
for any $i = 0,1,\dots$.
\medskip

The main point of this article is computational. We
realized that, while waiting for a more systematic exposition of the
general ideas, recalled above, we were risking to forget about some
difficult  computations we had previously made to support our  ideas.
We therefore decided to store those computations in these
proceedings, and at the same time to make them available to any
willing colleague.
\par


\bigskip
{\bf Contents.}
\vskip 5 pt
\item{1.} Filtered logarithmic $F$-crystals.
\item{2.} Dwork families of logarithmic $F$-crystals.
\item{3.} The hypergeometric family.
\itemitem{3.1} Determination of the Frobenius matrix.
\itemitem{3.2} The unit root subcrystal.
\itemitem{3.3} The Koblitz-Diamond formula.
\item{4.} Appendix: determination of Frobenius eigenvalues {\it via}
             modular relations.

\vskip 10 pt
{\bf \N {Filtered logarithmic $F$-crystals.}}
\vskip 5 pt

\par {\proc {\bf\NN.} {Notation. }}
Let $K$ be a complete discrete
valuation field of mixed characteristics $(0,p)$ with perfect residue
field
$k$, and let $\c V$ denote its ring of integers. We denote by
$| \;\;|$ the absolute value of $K$, normalized by $|p| = p^{-1}$.
For simplicity we will assume here that ${\cal V} = W(k)$ is
absolutely unramified, and we will denote by $\sigma: {\cal V} \to
{\cal V}$  the  Frobenius (and its extension to $K$).

\par {\proc
{\bf\NN.} {}}
\newbox\riferXD \global\setbox\riferXD=\copy\corrente
Let  $X$ be a smooth  $p$-adic formal scheme over  of finite type
over $\Spf \cal V$, and $D$ be
a divisor  in $X$ with strict normal
crossings relative to
$\cal V$.  The pair $(X,D)$ defines a fine log scheme $(X,M_D)$ in
the sense of Fontaine-Illusie and Kato  [KK], over  $\Spf \cal V$
endowed with the trivial log structure. It reduces modulo $p$ to a
similar pair $(X_k,D_k)$.
   In local \'etale coordinates $(x_1,
\dots,x_n)$, we may assume that the ideal sheaf ${\cal I}_D$ of $D$
is generated by $x_1 \cdots x_d$, and that the sheaf of log
differentials $\Omega_X^1(\log D) = \Omega_{(X,D)}^1$   admits the
${\cal O}_X$-basis
   ${{dx_1} \over {x_1}}, \dots, {{dx_d} \over
{x_d}}, dx_{d+1},\dots,dx_n$. We denote by ${\cal S}_D$ the ${\cal
O}_X$-algebra $\bigoplus_{j \geq 0}  {\cal I}^{-j}_D$.
For any
formal scheme $T$ over $\Spf \cal V$,
$T^{\sigma}$ will denote the
formal scheme over $\Spf \cal V$
obtained by the base change $\sigma$. Similarly,
for an ${\cal O}_T$-module ${\cal E}$ (with connection $\nabla$),
${\cal E}^{\sigma}$
will denote the ${\cal O}_{T^{\sigma}}$-module (with connection
$\nabla^{\sigma}$)
obtained by the base change $\sigma$.

\par {\proc {\bf\NN.} {Definitions.}}
A {\sl  logarithmic  crystal on $(X,D)/{\cal V}$ }
consists of
\item {$(a)$}
a finite projective ${\cal O}_X$-module $\cal E$;
\item {$(b)$}
an integrable logarithmic connection $\nabla:{\cal E} \to {\cal E}
\otimes_A\Omega^1_X(\log D)$.
\smallskip\noindent The logarithmic 
crystal $({\cal E}, \nabla)$ on $(X,D)/{\cal V}$ is 
{\sl convergent} 
if  the following convergence condition holds:
\item {$(c)$} for any $\epsilon > 0$, and any  section $e$ of $\cal E$ over a
coordinate domain $(U, x_1, \dots, x_n)$, the section of $\cal E$
defined by
$$
p^{\epsilon \Sigma_i\mu_i}{1\over {\Pi_i(\mu_i
!)}}\prod_{i=1}^n\prod_{j=0}^{\mu_i}(\nabla_i-j)e
$$
converges to $0$ when $\Sigma_i\mu_i$ goes to $\infty$.
Here  the operators $\nabla_i$ are defined by the formula
$\nabla(e)=\sum_i \nabla_i(e){{dx_i} \over x_i}$.
\smallskip\noindent
A (resp.  convergent) logarithmic crystal $({\cal E}, \nabla)$ on 
$(X,D)/{\cal V}$ is
{\it filtered} if it is equipped with
\item {$(d)$}
a decreasing and finite filtration, exhaustive and separated, $\{
\Fil^i {\cal E} \}_i$ by local direct factors of $\cal E$
satisfying
the Griffiths transversality
$\nabla(\Fil^i {\cal E}) \subseteq \Fil^{i{-}1} {\cal E} \otimes
\Omega^1_X (\log D)$.

\medskip
For $f \in {\b Z}_{\geq 1}$, a lifting of the $f$-th order relative
Frobenius of $X_k$ to the log scheme $(X,D)$   is
     a  lifting $\phi_f  : X \to X^{\sigma^f} $ of the $f$-th order
relative Frobenius of $X_k$  (coming from the $f$-th iterate of the
standard Frobenius)
{\it adapted} to the divisor $D$, {\it i.e. } such that
$\phi_f^{\ast} D^{\sigma^f} = p^fD$. So, in terms of local
coordinates as above, for any $i$ one has $\phi_f^{\ast}(x_i) =
x_i^{p^f}$ up to units in ${\cal O}_X$. When $f=1$, we avoid
mentioning $f$.
\par
A { \it logarithmic  $F$-crystal over
$(X,D)/{\cal V}$} for the $f$-th order Frobenius,  is a convergent logarithmic
crystal $({\cal E}, \nabla)$ over $(X,D)/{\cal V}$,  together with
the assignment,
for any local lifting of the $f$-th order Frobenius
$\phi_f: U \to U^{\sigma^f}$, adapted to $D \cap U$, to an open
formal $\cal V$-subscheme $U$ of $X$,  of a horizontal
${\cal
S}_{D\cap U}$-linear monomorphism
$$F(\phi_f): \phi_f^{\ast}( {\cal
S}_D\otimes {\cal E}, \nabla)^{\sigma^f}_{|U^{\sigma^f}} \to ( {\cal
S}_D \otimes {\cal E}, \nabla)_{|U}
$$
which becomes an isomorphism
when tensored with $K$. (Notice that $\nabla$ induces a logarithmic
connection on ${\cal S}_D \otimes {\cal E}$, and that $\phi_f$
extends to a morphism of  ringed spaces
$(U,  {\cal S}_{D\cap U} )
\to (U,  {\cal S}_{D\cap U} )^{\sigma^f}$.)
\par  We say that $({\cal E}, \{ \Fil^i {\cal
E}\}_i, \nabla, F)$ is a  {\it filtered logarithmic $F$-crystal
over $(X,D)/{\cal V}$}
for the $f$-th order Frobenius, if $({\cal E},  \nabla, F)$ is a 
logarithmic  $F$-crystal over $(X,D)/{\cal V}$
for the $f$-th order Frobenius, and  $({\cal E}, \{ \Fil^i {\cal
E}\}_i, \nabla)$ a filtered convergent logarithmic crystal over 
$(X,D)/{\cal V}$. Then, 
$({\cal E}, \{ \Fil^i {\cal
E}\}_i, \nabla, F)$
is {\it divisible}
if, for any $i$,
$$
F(\phi_f) \left(\phi_f^{\ast}  (\Fil^i {\cal
E}^{\sigma^f})_{|U^{\sigma^f}}\right) \subset p^{i} ({\cal S}_D \otimes
{\cal E})_{|U} \; .
$$
(In our application, the stronger condition $$
F(\phi_f) \left(\phi_f^{\ast}  (\Fil^i {\cal
E}^{\sigma^f})_{|U^{\sigma^f}}\right) \subset p^{fi} ({\cal S}_D \otimes
{\cal E})_{|U} \; ,
$$
will be considered, at the expense of the strength of our results.)
We omit the natural extensions of the previous
definitions to a relative situation $(X,D)/S$, where $S$ is a formal
$p$-adic scheme, smooth and of finite type over $\Spf {\cal V}$ and
$D$ is a divisor in $X$,  with strict normal crossings relative to
$S$. In the present discussion however a logarithmic crystal on 
$(X,D)/S$ will rarely be convergent in the natural relative 
generalization of this notion. This is because the points of $S$ will 
play for us the role of  variable ``exponents of monodromy" rather 
than the one of rational parameters as in the theory of Picard-Fuchs 
equations.

\smallskip
{\proc  {\bf\NN.} {}}
Let  $({\cal E}, \nabla, F)$ be  a  logarithmic
$F$-crystal  over $(X,D)/{\cal V}$ for the $f$-th order Frobenius,
and let us assume, for simplicity, that $X$ is of relative dimension
1 over $\cal V$, with coordinate $x$, and that $D$ consists of a
finite set of $\cal V$-valued points  $x = t_1, \dots, t_r$, with
$t_i \in {\cal V}$. In particular, as a topological space, $D$
consists of a finite set of $k$-valued points $\ov{x} = \ov{t}_1,
\dots, \ov{t}_r$, for the reduced coordinate $\ov{x}$. We recall the
{\it bounded Robba ring} ${\cal R}_{X, \ov{t}_i}$ of $X$ at
$\ov{t}_i$, which is the ring of Laurent  series  $\sum_{j \in {\b
Z}} a_j (x-t_i)^j$, with $a_j \in {\cal V}$, converging in some
annulus $\epsilon < |x-t_i| < 1$.  Then, according to [Ke, \S 4], one
can develop Dieudonn\'e theory
over ${\cal R}_{X, \ov{t}_i}$, and
consequently define the notion of {\it special Newton polygon} of
$({\cal E}, \nabla, F)$ at ${\ov t}_i$, that is over ${\cal R}_{X, {\ov t}_i}$.
So, combining this result with the classical theory of Newton
polygons of $F$-crystals [Ka], we associate to $({\cal E}, \nabla,
F)$  a Newton polygon at each $\ov{k}$-valued point of the
$k$-scheme $X_k$, where $\ov{k}$ denotes the algebraic closure of
$k$. We will say that $({\cal E}, \nabla, F)$  is a {\it unit root}
logarithmic
$F$-crystal  over $(X,D)/{\cal V}$, if its Newton polygon is a
horizontal segment at any $\ov{k}$-valued point of the
$k$-scheme $X_k$.

  \endgraf

\medskip
{\proc  {\bf\NN.} {}}
\newbox\riferKATZ \global\setbox\riferKATZ=\copy\corrente
We will give an example of the following reasonable generalization of
[Ka, 4.1], whose proof does
not seem to appear in the literature. We plan to give full details elsewhere.

\proclaim {Theorem}. {
Let $({\cal E}, \{\Fil^2{\cal E} = 0 \subset \Fil^1{\cal E} \subset
\Fil^0 {\cal E} = {\cal E}\}, \nabla, F)$
be a divisible $F$-crystal for
the $f$-th order Frobenius on $(X,D)/{\cal V}$, with a two-step
filtration. Assume $X$ is of relative dimension 1 over $\cal V$, and
that at every $\ov{k}$-valued point of $X_k$, its Newton polygon
begins with a side of slope zero, always of the same length $\nu \geq
1$ ({\it i.e.}, point by point, the unit root part has rank $\nu$)
and that ${\cal E}/\Fil^1 {\cal E}$ is of constant rank $\nu$. Then
there exists a logarithmic unit root $F$-sub-crystal $({\cal U}, \nabla,F)
\subset ({\cal E}, \nabla,F)$, whose underlying module $\cal U$ is
transversal to $\Fil^1{\cal E}$ ({\it i.e.} ${\cal E} = {\cal U}
\oplus \Fil^1{\cal E}$).
}
\par

\vskip 10 pt
\def\u{\vec}

{\bf \N {Dwork families of logarithmic $F$-crystals.}}
\vskip 5 pt
This section is an abstract formulation of the theory of generalized 
hypergeometric
functions of
[GHF] and [Ad].
\par
{\proc  {\bf\NN.} {}}
Let $(X,D)$ be as in the previous section, and let $H = \Spf {\b
Z}_p \{a_1, \dots , a_r \}$, $(X_H,D_H) = (X,D)
\times H$,
with projections $p_X: X_H \longrightarrow X$ and $p_H: X_H \to H$.
The coordinates $a_1, \dots , a_r$ will play a special role in what
follows, together with a finite
set of linear forms ${\cal L} = \{\ell_1({\u a}), \dots, \ell_N({\u
a}) \} \subset {\b Z}
[a_1,\dots,a_r]$.

We will assume that the system of  inequalities
$\ell_i({\u a}) \geq 0$, for $i =1, \dots, N$,  defines  a rational
polyhedral cone $C_{\cal L}$ of dimension $r$ in ${\b R}^r$, and
that, for any $i$, $\ell_i({\u a}) = 0$ is a 1-codimensional face of
$C_{\cal L}$,  and  $\ell_i ({\b Z}^r ) = {\b Z}$.
A {\it meromorphic function} on an open formal $\cal V$-subscheme of
$X_H$ will be assumed to have a finite set of polar hypersurfaces of
the form $p_X^{-1}(D)$ and of the form  $p_H^{-1}(\ell_i({\u a}) =
j)$, for $i = 1,\dots, N$ and $j \in {\b Z}$. So, in local
coordinates $(x_1, \dots, x_n)$ over an open formal subscheme $U$ of
$X$,  with $D \cap U = V(x_1 \cdots x_d)$, a meromorphic function $g$
on $U \times H$ will be a quotient of a section $h \in \Gamma(U
\times H, {\cal O}_{X_H})$, by an expression of the form $x_1^{u_1}
\cdots x_d^{u_d} \prod_{i,j} (\ell_i({\u a}) - j)$.
We will loosely
talk about meromorphic structures on $X_H$ in that sense.

\smallskip
{\proc  {\bf\NN.} {}}
Let ${\cal E}_0$ be a locally free ${\cal O}_X$-module
of finite type, ${\cal E} = p_X^{\ast}{\cal E}_0$, $\{ \Fil^i {\cal E}
\}_i$ be a filtration of ${\cal E}_0$ by local direct factors,
and let
$\nabla$ be an integrable
$(X_H, D_H)/H$-connection on $\cal E$,  such that  $({\cal E}, \{
\Fil^i{\cal E} \}_i, \nabla)$ becomes a filtered  logarithmic
crystal over $(X_H, D_H)/H$. We will replace the relative convergence 
condition by the following weaker convergence 
condition
\item{($c^{\prime}$)} For any ${\u a} \in H(\b Q \cap \b 
Z_p)$
$({\cal E}_{\u a},  \{(\Fil^i {\cal E})_{\u a}
\}_i , \nabla_{\u a}) =
{\u a}^{\ast}({\cal E}, \{\Fil^i {\cal E}
\}_i, \nabla)$ is a filtered convergent
logarithmic crystal.

\smallskip
{\proc  {\bf\NN.} {}}
For ${\u u} \in {\b G}_a^r({\b Z})$, we
denote by
$$\sigma_{\u u}:H \longrightarrow H
$$
the translation mapping ${\u a} \longmapsto {\u a}+{\u u}$, inducing
$$\sigma_{\u u} = \id_X \times \sigma_{\u u}:X_H \longrightarrow X_H \; .
$$
We assume the existence of horizontal meromorphic isomorphisms
$$b_{\u u} :  ({\cal E},   \nabla)
\longrightarrow \sigma_{\u u}^{\ast} ({\cal E}, \nabla)
$$
such that for ${\u u} , {\u v} \in {\b G}_a^r({\b Z})$
$$b_{{\u u}+{\u v}} = \sigma_{\u u}^{\ast}(b_{\u v}) \circ b_{\u u}\;.
$$
We will assume that for ${\u \mu} \in {\b Z}^r \cap C_{\cal L}$,
the only poles of the map $b_{\u u}$ be along $D$.  To be more
precise [GHF, Conjecture 6.3.1], [Ad, Thm. 8.1], we may consider the
natural map of ringed spaces $q_X: (X_H,
p_X^{\ast} {\cal S}_D) \to X_H$. Then, we will assume that, for ${\u
u} \in {\b Z}^r \cap C_{\cal L}$,  $b_{\u u}$ is a honest morphism
$$
b_{\u u} :  q_X^{\ast} ({\cal E},   \nabla)
\longrightarrow q_X^{\ast}  \sigma_{\u u}^{\ast} ({\cal E}, \nabla)
\; ,
$$
and  that $\det b_{\u u}$  vanishes  if and only if
$$\prod_{i=1}^N(\ell_i({\u
a}))_{\ell_i({\u u})} = 0 \; .
$$

\par

{\proc  {\bf\NN.} {}}
For
${\u \mu} \in {\b Z}^r$ and $\rho = p^{-s}$, $s \in {\b Z}_{\geq 0}$,
we define $\b D(-{\u \mu}, \rho)$ to be the formal counterpart of the
closed analytic disk $D(-{\u \mu}, \rho) = \{ {\u a} \in {\b C}_p^r \mid
|a_i +\mu_i| \leq \rho, \forall i = 1,\dots,r \}$, that is the formal
$\cal V$-subscheme of $H$
$$\b D(-{\u \mu}, \rho) = \Spf {\cal V} \{
a_1 ,\dots,   a_r, b_1 ,\dots,  b_r  \} / (p^s b_1 - a_1 - \mu_1,
\dots, p^s b_r - a_r - \mu_r) \; .$$
For $\rho < 1$, there are
natural morphisms
$$\eqalign{
{\tau_{\u \mu}}: \b D(-{\u \mu}, \rho) &\longrightarrow
\b D({\u 0}, p\rho) \cr
{\u a} & \longmapsto {{{\u a}+{\u \mu}}\over p} \; .
\cr}
$$
Notice that for any ${\u \mu}, {\u u}, {\u v} \in {\b G}_a^r({\b Z})$
$$\sigma_{\u u} \circ \tau_{\u \mu} =
\tau_{p{\u u}-{\u v}+{\u \mu}} \circ \sigma_{\u v} \; .
$$

\par
Let $ U$ be an  open formal subscheme  of $X$ and
$\phi : U \longrightarrow U^{\sigma}$ be a
$\sigma$-linear lifting
of Frobenius adapted to $U \cap D$. We consider the map
$$\eqalign{
\phi \times \tau_{\u \mu} : U \times \b D (-{\u \mu},\rho)
&\longrightarrow U^{\sigma} \times \b D ({\u 0},p\rho)  \cr
(x,{\u a}) &\longmapsto \left(\phi(x), {{{\u a}+{\u \mu}}\over p}\right) \;.
\cr}
$$
We assume that for each $(U,\phi)$, $\rho \in p^{{\b Z}_{< 0}}$, and
${\u \mu} \in {\b Z}^r$, as before,
there exists a
{\it meromorphic} morphism of  logarithmic crystals over
$(U \times {\b D (-{\u \mu},\rho)}, D \times {\b D (-{\u 
\mu},\rho)})/\b D (-{\u
\mu},\rho)$
$$
F({\phi},{\u \mu}) : (\phi \times \tau_{\u \mu})^{\ast}
({\cal E}^{\sigma})_{|U\times \b D (-{\u \mu},\rho)}
\longrightarrow {{\cal E}}_{|U \times \b D (-{\u \mu},\rho)} \; .
$$
The previous data should satisfy
$$
\sigma_{\u v}^{\ast}(F(\phi,p{\u
u}-{\u v}+{\u \mu})) \circ (\phi \times \tau_{\u \mu})^{\ast}(b_{\u
u}^{\sigma}) = b_{\u v} \circ F(\phi,{\u \mu}) \; .
$$
  What this simply means is that  the map $F({\phi},{\u \mu})$ may
be coherently regarded as a system of maps
  $$F({\phi},{\u \mu})_{(x,
{\u a})} = F({\u a}, {\u b}; x, \phi (x)): %
{\cal E}^{\sigma}_{({\u b}, \phi (x))} \to
{\cal E}_{({\u a}, x)}  \; , $$
  for any, say,  $\cal V$-valued points
$(x, {\u a})$ of $U \times \b D (- {\u \mu}, \rho)$ and $(\phi (x), 
{\u b})$  of $U^{\sigma} \times \b D ({\u 0}, p \rho)$,  with $p{\u 
b} - {\u a} = {\u \mu} \in {\b Z}^r$.
  It will also be
sometimes convenient to use the notation $F(\phi ; {\u a}, {\u b}; 
x)$ or $F(\phi ; {\u a}, {\u b})$ for that map. 
From the viewpoint 
of $p$-adic convergence, our condition means that the matrix 
$\gamma^{(\phi)}({\u a}, {\u b}; x)$ expressing 
$F(\phi ; {\u a}, 
{\u b}; x)$ in terms of a global basis of ${\cal E}_0$, is 
$p$-adically meromorphic in the variables $({\u a}, {\u b}; x)$, for 
{\it fixed } ${\u \mu} = p{\u b} - {\u a}$  in $\b Z^r$ and ${\u b} 
\in D (\u 0, 1)$.

\par
A more precise meromorphy condition, satisfied for generalized
hypergeometric functions [GHF, 4.7.1], [Ad, 9.12],  is that the polar locus
of $F({\phi},{\u \mu})$ should be contained in the union of $p_X^{-1} 
(D)$ and of the zero locus of the
determinant of the map $(\phi \times \tau_{\u \mu})^{\ast}(b_{\u
v})$, where ${\u v} \in {\b Z}^r \cap C_{\cal L}$ is such that
$${\b
Z}^r \cap ({\u v} - {  {\u \mu} \over p} + C_{\cal L} ) \subset
C_{\cal L} \; .
$$
In particular, [GHF, 6.13.2], [Ad, 9.14],  if  ${\u \mu}
\in {\b Z}^r \cap C_{\cal L}$, and $\ell_i({\u \mu}) \leq p-1$, for
any $i =1, \dots, N$, then $F({\phi},{\u \mu}) $ is a honest morphism
$$
F({\phi},{\u \mu}) : ((\phi \times \tau_{\u \mu})^{\ast}
(q_X^{\ast}{\cal E})^{\sigma})_{|U\times \b D (-{\u \mu},\rho)}
\longrightarrow ({q_X^{\ast}{\cal E}})_{|U \times \b D (-{\u \mu},\rho)} \; .
$$

\smallskip
{\proc  {\bf\NN.} {Definition.}}
A set of data $({\cal E}, \{\Fil^i {\cal E} \}_i,\nabla, \{
b_{\u u} \}_{\u u}, F)$ as before,
will be called
a {\it Dwork family of filtered convergent logarithmic $F$-crystals
on $X$, parametrized by $H$,
with set of singular forms $\cal L$,
on $(X,D)/{\cal V}$}.
  \par

\smallskip
{\proc  {\bf\NN.} {}}
  The notion of divisibility for Dwork family of filtered
$F$-crystals  $({\cal E}, \{\Fil^i {\cal E} \}_i,\nabla, \{
b_{\u u} \}_{\u u}, F)$, is perhaps a little unexpected. It first
requires that the filtration $\{\Fil^i {\cal E} \}_i$ is constructed
in the following way. Let  ${\cal A}= \{ {\u \mu} \}$ be a set of
representatives of $H({\b Z}_p)$ modulo $p$. Therefore $H({\b Z}_p)$
is the disjoint union of $\b D (-{\u \mu}, p^{-1})$, for ${\u \mu}  \in
{\cal A}$. We will assume to be  given a family, indexed by ${\u \mu}
\in {\cal A}$, of filtrations $\{\Fil^i_{\u \mu} {\cal E}_0 \}_i$  of
${\cal E}_0$ by local direct factors. We assume that, for all $i$ and
$\u \mu$,  on $X \times \b D (-{\u \mu}, p^{-1})$, $\Fil^i {\cal E} $
coincides with the inverse image via the first projection of
$\Fil^i_{\u \mu} {\cal E}_0$.
  \par
The Dwork family will then said
to be  {\it divisible} if, for any $i$,  any ${\u \mu} \in {\b Z}^r$,
and any $\rho = p^{-s} < 1$,
  $$
F({\phi},{\u \mu}) \left(((\phi \times \tau_{\u \mu})^{\ast}
(\Fil^i {\cal E})^{\sigma})_{|U\times \b D (-{\u \mu},\rho)}\right)
\subset p^i({q_X^{\ast}{\cal E}})_{|U \times \b D (-{\u \mu},\rho)} \; .
$$
\par

\smallskip
{\proc  {\bf\NN.} {}}
Let $({\cal E}, \{\Fil^i {\cal E}
\}_i,\nabla, \{ b_{\u u} \}_{\u u}, F)$ be a Dwork family of 
logarithmic convergent filtered
$F$-crystals on $(X,D)/ {\cal V}$, parametrized by $H = \Spf {\b
Z}\{a_1, \dots, a_r \}$, with set of singular forms ${\cal L} = \{
\ell_1, \dots, \ell_N \}$.
\endgraf
Now, for  ${\u a} \in {\b Z}_p^r$, let us choose
two  sequences ${\u a}^{[i]} \in {\b Z}_p^r$ and ${\u \mu}^{[i]}  \in
{\b Z}^r$, $i =0,1,\dots$, so that ${\u a}^{[0]}  = {\u a}$,
   and
$p{\u a}^{[i+1]}
- {\u a}^{[i]}  = {\u \mu}^{[i]} $, for any $i$. One possible choice
is ${\u a}^{[i]}  ={\u a}^{(i)} = (a_1^{(i)}, \dots ,a_r^{(i)})$ and 
${\u \mu}^{[i]}  = {\u \mu}_{{\u
a}^{[i]} } = (\mu_{a_1^{(i)}}, \dots, \mu_{a_r^{(i)}})$, as defined 
in the introduction.
If ${\u a} \in (\b Q \cap {\b Z}_p)^r$, then 
there exists
$f \in {\b Z}_{\geq 1}$  such that $(p^f-1){\u a} \in {\b Z}^r$:
 the minimal 
such $f$ will be called the {\sl period} of $\u a$, and we will say 
that $\u a$ is {\sl of finite period} $f$. So, if $\u a$ is  of 
finite period $f$, we may
arrange the previous choices so that
   ${\u a}^{[f]} = {\u a}$.

If $U$ is an open formal subscheme of $X$ and $\phi$ is a lifting of 
Frobenius on $U$, adapted to $D \cap U$, the map
$$
F(\phi,{\u \mu}^{[i]} ): (\phi \times \tau_{{\u \mu}^{[i]} })^{\ast}
({\cal E}^{\sigma})_{|U \times \b D (-{\u \mu}^{[i]} ,\rho)}
\longrightarrow
{\cal E}_{| U \times \b D (-{\u \mu}^{[i]} ,\rho)}
$$
is represented, in terms of
a basis of global sections $\u e$ of ${\cal E}_0$ over $U$, by a
matrix of functions meromorphic
in $U \times \b D (-{\u \mu}^{[i]} ,\rho)$. So, for ${\u a} \in 
H({\cal V})$, outside of a well-understood polar locus ({\it e.g.} if 
$\ell_i ({\u a}) \notin \b Z$, for $i =1, \dots,N$), the previous map
can be specialized to
induce a meromorphic map,
necessarily horizontal,
$$
F(\phi; {\u a}^{[i]} ,{\u a}^{[i+1]} ):(\phi^{\ast}
({\cal E}_{{\u a}^{[i+1]} }, \nabla_{{\u a}^{[i+1]} })^{\sigma})_{|U}
\longrightarrow ({\cal E}_{{\u a}^{[i]} }, \nabla_{{\u a}^{[i]} })_{|U}\;.
$$
If moreover ${\u a}$ is in $H(\b Q \cap \b Z_p)$, say ${\u a}^{[f]} = 
{\u a}$, our assumptions imply that
$$({\cal E}_{{\u a}^{[f]}}, \nabla_{{\u a}^{[f]}})
= ({\cal E}_{\u a}, \nabla_{\u a})\;,$$
and,  outside of some polar locus ({\it e.g.} if $\ell_i ({\u 
a}^{[j]}) \notin \b Z$, for $i =1, \dots,N$ and $j =0, \dots,f-1$), 
that
$$
F(\phi; {\u a}^{[f-1]},{\u a}^{[f]}) \circ
F(\phi;{\u a}^{[f-2]},{\u a}^{[f-1]})\circ \cdots
\cdots \circ
F(\phi;{\u a}^{[0]},{\u a}^{[1]})
$$
is a meromorphic horizontal map
$$F(\phi^f;{\u a}):(\phi^f)^{\ast}({\cal E}_{{\u a}}, \nabla_{{\u
a}})^{\sigma^f}_{|U}
\longrightarrow
    ({\cal E}_{\u a}, \nabla_{\u a})_{|U} .
$$
The structure $({\cal E}_{\u a}, \nabla_{\u a},F( \; - \;  ;{\u a}))$
is then a logarithmic  $F$-crystal on $(X,D)/{\cal V}$
for the $f$-th iterate of the Frobenius, in the usual sense.

The convergence assumption now implies, by [Ka, 3.1.2],  that  the
solutions of
$({\cal E}_{\u a}, \nabla_{\u a})$
    at any rigid point  $x$ of the Raynaud generic fiber
$X_K$ of $X$, not in the open tube of $D_k$, converge in the open
tube of radius 1 around $x$. Therefore, by
[BC1],  the holomorphic part of the solution matrix at a point of
$D_K$, in the sense of the classical theory of regular singularities,
converges in the open tube of $D_k$.

\par

\medskip
{\proc  {\bf\NN.} {}}
\newbox\riferKATZGEN \global\setbox\riferKATZGEN=\copy\corrente
We now propose an extension of  \copy\riferKATZ \  for divisible Dwork
families of
logarithmic $F$-crystals. This paper gives an example of this 
situation in relative dimension 1.
\proclaim {{Question}}.
{Let $(X,D)$ be as in \copy\riferXD. Let $({\cal E}, \{\Fil^2{\cal E} 
= 0 \subset \Fil^1{\cal E} \subset
\Fil^0 {\cal E} = {\cal E}\}, \nabla, \{b_{\u u} \}_{{\u u} \in {\b
Z}^r}, F)$ be a divisible Dwork family of  logarithmic convergent filtered
$F$-crystals on $(X,D)/ {\cal V}$, parametrized by $H = \Spf {\b
Z}\{a_1, \dots, a_r \}$, with set of singular forms ${\cal L} = \{
\ell_1, \dots, \ell_N \}$.
  Let us assume that, for any
  ${\u a} \in H({\b Z}_p)$ of finite period $f$,  such that 
$\ell_i({\u a}) \notin \b Z$, for $i = 1, \dots, N$, the divisible 
logarithmic $F$-crystal
$({\cal E}_{\u a}, \{\Fil^2{\cal E}_{\u a} = 0 \subset \Fil^1{\cal E}_{\u a}
\subset
\Fil^0 {\cal E}_{\u a} = {\cal E}_{\u a}\}, \nabla_{\u a},  F(\; - \; ,{\u
a}))$ over $(X,D)/{\cal V}$ for the $f$-th order
Frobenius,   admits a
logarithmic unit root $F$-sub-crystal $({\cal U}_{\u a}, \nabla_{\u a},F(\; -
\; ,{\u a}))$, whose underlying module ${\cal U}_{\u a}$ is
transversal to
$\Fil^1{\cal E}_{\u a}$ ({\it i.e.} ${\cal E}_{\u a} = {\cal U}_{\u a}
\oplus \Fil^1{\cal E}_{\u a}$).
Then does there exist a
logarithmic sub-crystal $({\cal U}, \nabla,F)$ of $({\cal E},
\nabla)$ on $(X_H,D_H)/H$, stable under the map $F$, in the sense 
that 
 for any $(U,\phi)$ as above, for
any ${\u \mu} \in {\b Z}^r$,
and any $\rho = p^{-s} < 1$, $F$ induces 
a meromorphic morphism of 
logarithmic crystals over
$(U \times {\b D (-{\u \mu},\rho)}, D \times {\b D (-{\u 
\mu},\rho)})/\b D (-{\u
\mu},\rho)$
$$
F({\phi},{\u \mu}) : (\phi \times \tau_{\u \mu})^{\ast}
({\cal U}^{\sigma})_{|U\times \b D (-{\u \mu},\rho)}
\longrightarrow {{\cal U}}_{|U \times \b D (-{\u \mu},\rho)} \;  ,
$$
 whose underlying module $\cal U$
is transversal to
$\Fil^1{\cal E}$ ({\it i.e.} ${\cal E} = {\cal U} \oplus \Fil^1{\cal
E}$) and whose specialization at any  ${\u a} \in H({\b Q} \cap {\b 
Z}_p)$ coincides with
$({\cal U}_{\u a}, \nabla_{\u a},F(\; - \; ,{\u a}))$? }
\par

We remark that the rank of the underlying modules
${\cal U}_{\u a}$ may vary with the class of $\u a$ mod $p$. We treat 
a simplified version of the previous question, where we restrict our 
parameters ${\u a}$ to a bizarre subset  ${\cal T}_2$ of $H(\b Z_p)$, 
stable under the map ${\u a} \longmapsto {\u a}^{(1)}$, and such that 
the filtration $\Fil^2{\cal E}_{\u a} = 0 \subset \Fil^1{\cal E}_{\u 
a}
\subset
\Fil^0 {\cal E}_{\u a} =  {\cal E}_{\u a} $ of $ {\cal E}_{\u a}$ is 
independent of  ${\u a} \in {\cal T}_2$.
\vskip 10 pt
{\bf \N{The hypergeometric family. }}
\vskip 5 pt
A very interesting example of the preceding situation is connected 
with the hypergeometric
system:
$$
    {{dY} \over {d\lambda}} = Y G_{\vec a}(\lambda ),
     \qquad\qquad \vec a = (a_1,a_2,a_3) \in \b Z_p^3
$$
where
$$
    G_{\vec a}(\lambda )=
    \pmatrix{ -\ds{a_3 \over \lambda} & \ds{{a_3 - a_1} \over {1-\lambda}} \cr
              \ds{{a_3 - a_2} \over {\lambda}}
               & \ds{{a_1 + a_2 - a_3} \over {1-\lambda}} \cr}\ .
$$
In this case, the relevant linear forms are 
 $${\c L} = \{ 
\ell_1({\u a}) = a_3 - a_1  \;  , \;  \;  \ell_2({\u a}) = a_3 - a_2 
\;  , \;  \; \ell_3({\u a}) = a_2  \;  , \;  \;  \ell_4({\u a}) = a_1 
\} \; .$$
 We denote by $C_{\vec a}(z, \lambda )$ the matrix solution 
at
$z\not = 0,1, \infty$ such that $C_{\vec a}(z,z)= \b I_2$, the $2 
\times 2$ identity matrix.
When the entries of $\u a$ are rational, that matrix  converges if
$|\lambda - z | < |z| \min \l( 1, |1-z| \r)$ [LDE].

{\proc  {\bf\NNN.} {}}
\newbox\riferB \global\setbox\riferB=\copy\corrente
The $F$-crystal structure of the hypergeometric system is expressed 
by the following data.
We take $\vec a, \vec b \in \b Z_p^3$,
$p \vec b - \vec a = \vec \mu \in \b Z^3$,
$\phi$ any lifting of Frobenius on an open formal subscheme $U$ of 
$\hat {\b P}$, adapted to $ D = \{ 0,1, \infty \}$. We obtain a 
matrix
$\gamma^{(\phi)} (\vec a, \vec b ; \lambda)$ meromorphic in the 
variables $({\u a}, {\u b},\lambda)$, {\it fpr fixed $\u \mu$},  such 
that:
$$
    C_{\vec b}^\sigma \l( \phi (z), \phi (\lambda) \r)
    \gamma^{(\phi)} (\vec a, \vec b ; \lambda)^t =
    \gamma^{(\phi)} (\vec a, \vec b ; z)^t
    C_{\vec a}(z, \lambda ) \; ,
$$
whenever $C_{\vec a}(z, \lambda )$ converges; this holds, in 
particular,  if ${\u a} \in (\b Q \cap \b Z_p)^3$, $z \not = 0,1, 
\infty$ and
$|\lambda - z | < |z| \min \l( 1, |1-z| \r)$.
Notice that we transpose the Frobenius matrix
$\gamma^{(\phi)} (\vec a, \vec b ; \lambda)$
in this formula in view of compatibility with the notation of
Dwork in [Ku] and [GHF].

{\proc  {\bf\NNN.} {}}
\newbox\riferA \global\setbox\riferA=\copy\corrente
We assume in our calculations that the local monodromy semisimple. 
 
In  the singular classes we
formally write the solutions in the form:

$$\matrix{
    \pmatrix {1 &0 \cr 0 & \lambda^{-a_3} \cr} U_{\vec a}^{(0)}(\lambda )
    \hfill & \qquad\quad
    U_{\vec a}^{(0)}(0) = \pmatrix{ a_3 - a_2 & a_3 \cr 1-a_3 & 0 \cr }
    \hfill\cr
    \pmatrix {1 &0 \cr 0 & (1-\lambda)^{a_3 - a_1 - a_2} \cr}
    U_{\vec a}^{(1)}(\lambda )
    \hfill &\qquad\quad
    U_{\vec a}^{(1)}(1) =
    \pmatrix{ a_1 + a_2 - a_3 & a_1 - a_3 \cr 0 & a_1 + a_2 -a_3 - 1 \cr }
    \hfill\cr
    \pmatrix {\lambda^{-a_1} & 0 \cr 0 & \lambda^{-a_2} \cr}
    U_{\vec a}^{(\infty)}(\lambda )
    \hfill &\qquad\quad
    U_{\vec a}^{(\infty)}(\infty) =
    \pmatrix{ a_3 - a_2 & a_3 - a_1 \cr a_2 - a_1 + 1 & a_2 - a_1 + 1 
\cr } \; \; .
    \hfill\cr}
$$
Under suitable assumptions on ${\u a}$, $U_{\vec a}^{(i)}(\lambda )$ 
is a holomorphic matrix on the residue
class of $i \in \{ 0,1, \infty \}$.
Let
$U_{\vec a}^{(i)}(\lambda ) =
{
   \pmatrix{u_1^{(i)} & u_2^{(i)}\cr u_3^{(i)} & u_4^{(i)}\cr }}$
for $i \in \{ 0,1, \infty \}$.
We know that
$$\hskip-40pt
\l\{ \matrix{ u_1^{(0)} = (a_3 -a_2)F(a_1,a_2,a_3 +1;\lambda) \hfill\cr
                  u_2^{(0)} = a_3 F(a_1,a_2,a_3 ;\lambda) \hfill\cr
                  u_3^{(0)} = (1-a_3)
                    F(a_2-a_3, a_1-a_3, 1-a_3; \lambda) \hfill\cr
                  u_4^{(0)} = (a_3 - a_1)\lambda
                    F(1+a_2 - a_3, 1+a_1 - a_3, 2-a_3; \lambda) \hfill\cr}
    \r.
$$
$$ \hskip 40pt
\l\{ \matrix{ u_1^{(1)} = (a_1+a_2 -a_3)F(a_1,a_2,a_1+a_2 -a_3;1-\lambda)
                  \hfill\cr
                  u_2^{(1)} = (a_1 - a_3) F(a_1,a_2,a_1+a_2 -a_3 +1 ;1-\lambda)
                  \hfill\cr
                  u_3^{(1)} = (a_3-a_2)(1-\lambda)
                   F(a_3-a_1+1, a_3-a_2+1,a_3-a_1-a_2+2;1-\lambda)
                  \hfill\cr
                  u_4^{(1)} = (a_1 +a_2-a_3-1)
                  F(a_3-a_1,a_3-a_2, a_3-a_1-a_2+1;1- \lambda)
                 \hfill\cr}
    \r.
$$
$$\hskip-40pt
\l\{ \matrix{ u_1^{(\infty)}
                  = (a_3 -a_2)F(a_1,a_1-a_3,a_1-a_2 +1;\lambda^{-1}) \hfill\cr
                  u_2^{(\infty)}
                  = (a_3-a_1) F(a_1,a_1-a_3+1,a_1-a_2+1 ;\lambda^{-1}) \hfill\cr
                  u_3^{(\infty)}
                  = (a_2-a_1+1)
                  F(a_2-a_3, a_2, a_2-a_1+1; \lambda^{-1}) \hfill\cr
                  u_4^{(\infty)}
                  = (a_2 - a_1+1)
                  F(a_2 - a_3+1,a_2,a_2-a_1+1; \lambda^{-1})
                  \hfill\cr}
    \r.
$$
({\it cf.} [Po, \S 22],  or, more precisely, [LDE, Lemma 24.1] for 
the solution at $0$).

{\proc {\bf\NNN.} {}}
We denote by  $\phi_i$ a lifting of Frobenius to a formal 
neighborhood $U_i$ of  $i \in \{ 0,1, \infty \}$, adapted to $i$.
Since  
 $\phi_0(\lambda ) = \lambda^p \l( 1 + \epsilon\lambda u(\lambda ) \r)$
with $|\epsilon | < 1$ and $u(\lambda ) \in \c V \aqq \lambda \cqq $,
we have that
$$ {{\phi_0(\lambda )^{b_3}} \over {\lambda^{a_3}}} =
     \lambda^{\mu_3} \l( 1 + \epsilon\lambda v(\lambda ) \r)
$$
where $v(\lambda ) \in \c V \aqq \lambda \cqq $
converges inside a disk of radius $>1$, and
  $$
    U_{\vec b}^{(0)\,\sigma} \l( \phi_0(\lambda ) \r)
    \gamma^{(\phi_0)} (\vec a, \vec b ; \lambda)^t =
    \pmatrix{ \xi_1^{(0)} (\vec a, \vec b ) & 0 \cr
              0 & \xi_2^{(0)} (\vec a, \vec b )
                  \ds{{\phi_0(\lambda )^{b_3}} \over {\lambda^{a_3}}} \cr}
    U_{\vec a}^{(0)}(\lambda ) \; .
$$

 Similarly,
$$
    U_{\vec b}^{(1)\,\sigma} \l( \phi_1(\lambda ) \r)
    \gamma^{(\phi_1)} (\vec a, \vec b ; \lambda)^t =
    \pmatrix{ \xi_1^{(1)} (\vec a, \vec b ) & 0 \cr
              0 & \xi_2^{(1)} (\vec a, \vec b )
                \ds{{\l( 1 - \phi_1(\lambda ) \r)^{b_1 + b_2 - b_3}}
                    \over {(1-\lambda )^{a_1 + a_2 - a_3}}} \cr}
    U_{\vec a}^{(1)}(\lambda )
$$
and
$$
    U_{\vec b}^{(\infty)\,\sigma} \l( \phi_\infty (\lambda ) \r)
    \gamma^{(\phi_\infty )} (\vec a, \vec b ; \lambda)^t =
    \pmatrix{ \xi_1^{(\infty )} (\vec a, \vec b )
               \ds{{\lambda^{a_1}} \over {\phi_\infty (\lambda )^{b_1}}}& 0 \cr
              0 & \xi_2^{(\infty )} (\vec a, \vec b )
                  \ds{{\lambda^{a_2}} \over {\phi_\infty (\lambda )^{b_2}}} \cr}
    U_{\vec a}^{(\infty )}(\lambda )\ .
$$

Inspection of the dominant terms at $\lambda = 0,1,\infty$,  makes 
it clear that, for any $i = 0,1,\infty$ and $j =1,2$, the functions 
$\xi^{(i)}_j(\vec a, \vec b )$ have the same  $p$-adic meromorphic 
behavior as the entries of the function $\gamma^{(\phi_i)} (\vec a, 
\vec b ; \lambda)$, with possibly some extra poles. This might be 
surprising, since the matrix $ U_{\vec a}^{(i)}(\lambda )$ itself is 
{\it not} a $p$-adic meromorphic function of $({\u a}, \lambda)$.

\smallskip
{\bf \NN. {Determination of Frobenius matrix.}}

{\proc {\bf\NNN.} {}} 
For the calculations we use  notation and  results of [Ku].
In that paper, Dwork investigates the effect of Kummer transformations
on the solutions of the hypergeometric equation.
We need a generalized form of the theorem in [Ku, \S 4], where  only the
standard lifting of Frobenius  $\lambda \mapsto \lambda^p$ is
considered. We rewrite the original statement as
$$ H_m( \vec a , \lambda )
     \gamma \l( M_m(\vec a), M_m(\vec b) ;
      \theta_m(\lambda),\theta_m(\lambda^p)  \r) =
     \gamma (\vec a , \vec b , \lambda ) H_m^\sigma( \vec b , \lambda^p )
$$
(this is possible once we make explicit
$h_m(\vec a , \vec b ,\lambda) = h_m( \vec a , \lambda ) \big/
   h_m^\sigma( \vec b , \lambda^p )$
   and $H_m( \vec a , \lambda )=h_m( \vec a , \lambda )N_m(\lambda)$
in the original formula
$$ h_m( \vec a , \vec b ,\lambda )
     \gamma \l( M_m(\vec a), M_m(\vec b) ;
      \theta_m(\lambda),\theta_m(\lambda^p)  \r) = N_m(\lambda)^{-1}
     \gamma (\vec a , \vec b , \lambda )  N_m(\lambda^p)
$$
of Dwork's article).
Our generalized statement is then
\proclaim{Theorem}. {
Let $\lambda \mapsto \phi (\lambda)$ be a function analytic
in a region of the type
$$ \inf \l\{ |\lambda |, |\lambda^{-1}|, |1-\lambda | \r\} > 1-\epsilon $$
for some $\epsilon>0$ and ``close to Frobenius'', in the sense  that
$$\l| \phi (\lambda) - \lambda^p \r| < 1$$
for every $\lambda$ for which $|\lambda | = 1= |1-\lambda |$.
Then
$$ H_m( \vec a , \lambda )
     \gamma \l( M_m(\vec a), M_m(\vec b) ;
      \theta_m(\lambda),\theta_m(\phi (\lambda))  \r) =
     \gamma (\vec a , \vec b ; \lambda , \phi (\lambda))
      H_m^\sigma( \vec b , \phi (\lambda) )
$$
}\par
 
We also  write the previous formula in the form
$$ H_m( \vec a , \lambda )
     \gamma^{(\phi)} \l( M_m(\vec a), M_m(\vec b) ;
      \theta_m(\lambda)  \r) =
     \gamma^{(\phi)} (\vec a , \vec b ; \lambda )
      H_m^\sigma( \vec b , \phi (\lambda) ) \ .
$$

{\proc{Proof. }} We  first use formula (3.1.7) of [Ku] for the
variation of the lifting of Frobenius, we then apply the original 
statement, then
  formula (2.9) of [Ku] (behaviour of solutions under Kummer transformations)
and finally again (3.1.7) of [Ku]:
$$
\eqalign{ &\quad
    H_m( \vec a , \lambda )
     \gamma^{(\phi)} \l( M_m(\vec a), M_m(\vec b) ; \theta_m(\lambda)  \r) = \cr
    &= H_m( \vec a , \lambda )
       \gamma^{(\lambda \mapsto \lambda^p)}
         \l( M_m(\vec a), M_m(\vec b) ; \theta_m(\lambda)  \r)
         C_{M_m(\vec b)}
         \l( \theta_m(\phi (\lambda)),\theta_m(\lambda^p) \r)^t =\cr
    &= \gamma^{(\lambda \mapsto \lambda^p)} (\vec a , \vec b ; \lambda )
        H_m^\sigma( \vec b , \lambda^p )
        C_{M_m(\vec b)}
         \l( \theta_m(\phi (\lambda)),\theta_m(\lambda^p) \r)^t =\cr
    &= \gamma^{(\lambda \mapsto \lambda^p)} (\vec a , \vec b ; \lambda )
        C_{\vec b}(\phi (\lambda) , \lambda^p )^t
        H_m^\sigma( \vec b , \phi (\lambda) ) = \cr
    &= \gamma^{(\phi)} (\vec a , \vec b ; \lambda )
        H_m^\sigma( \vec b , \phi (\lambda) ) \ .\cr
}$$
\hfill$\square$

\medskip
\par {\proc {\bf\NNN.} {Notation. }}
In the following formulas we use  Dwork's  symbol
$\gp (x,y)$ defined in [LDE, ch. 21],
for $py-x=\mu \in \b Z$ and $d(x,\b Z) \leq p^{-1}$, by
$$
    \gp (x+m,y+n)= (-\pi)^{n-m}
     {{\G (x+m) \G(y)}\over {\G(x) \G(y+n)}} \gp (x,y)
$$
for any $m,n \in \b Z$ and
$$ \gp (x,y)  = \pi^\mu \Gp (x) $$
if $\mu \in \{ 0, 1, \dots , p-1 \}$.
Here $\pi$ is a fixed element in $\o{\b Q}_p$ such that $\pi^{p-1} = -p$.
We recall also the symplectic relation
$$ \gp (x,y) \gp (1 - x,1 - y) = \m^{py-x}p $$
(equivalent to $ \Gp (x) \Gp (1 - x) = -\m^t$ with
$t\equiv -x \ \mod p$, $t \in \{ 0,\dots , p-1\}$).

\bigskip

{\proc {\bf\NNN.} {}}  For the Frobenius matrix at the origin,
Dwork [LDE, ch. 25] obtains the values:
$$ \eqalign{
    \xi_1^{(0)} (\vec a, \vec b ) &=
    {{\gp (a_2,b_2) \gp (a_3 - a_2, b_3 - b_2) }
      \over
     {\gp ( 1 + a_3, 1+ b_3 ) }} \cr
    \xi_2^{(0)} (\vec a, \vec b ) &= \m^{\mu_2 - \mu_3}
    {{\gp (a_3 - 1, b_3 - 1 ) \gp (1- a_1, 1- b_1) }
      \over
     {\gp (1 + a_3 - a_1, 1 + b_3 - b_1 ) }} \;  .\cr }
$$
We complete his calculations:

\proclaim {{\bf\NNN.} {Theorem}}. {
$$ \eqalign{
    \xi_j^{(1)} (\vec a, \vec b ) &=
     \m^{\mu_2} \xi_j^{(0)} \l( M^{(1)}(\vec a), M^{(1)}(\vec b) \r) \cr
    \xi_j^{(\infty)} (\vec a, \vec b ) &= \m^{\mu_1 + \mu_2 - \mu_3}
      \xi_j^{(0)} \l( M^{(\infty)} (\vec a), M^{(\infty)} (\vec b) \r) \cr }
$$
where $j=1,2$ and
$$ \eqalign{
         M^{(1)}(\vec a) &= (a_1,a_2, a_1 + a_2 - a_3 ) \cr
         M^{(\infty)} (\vec a) &= (a_1, a_1 - a_3, a_1 - a_2 ) \ .\cr  }
$$}\par
{\proc Proof.}
The argument is to compare the Frobenius action at infinity with
the Frobenius action at the origin by using the transformation $\theta_9$
of [Ku];
then to compare the Frobenius at $1$ with the one  at infinity by
using the transformation $\theta_{11}$ of [Ku].
We need for a lemma relating the analytic part of solutions
at a point and its image:
\proclaim {{\bf\NNN.} {Lemma}}. {
$$ \eqalign{
      U_{\vec a}^{(\infty)}(\lambda )
      &=U_{M_{9}({\vec a})}^{(0)}(\theta_9(\lambda) ) N_9^t\cr
      U_{\vec a}^{(1)}(\lambda )
      &=U_{M_{5}({\vec a})}^{(0)}(\theta_5(\lambda) ) N_5^t\cr
      U_{\vec a}^{(1)}(\lambda )
      &=U_{M_{11}({\vec a})}^{(\infty)}(\theta_{11}(\lambda) ) N_{11}^t  \cr
   }$$
where the transformations are named after  Dwork's tables in [Ku].
In particular $\theta_9(\lambda) = \lambda^{-1}$,
$\theta_5(\lambda) = 1 - \lambda$,
$\theta_{11}(\lambda) =( 1 - \lambda)^{-1}$.
}\par
{\proc {Proof (Lemma). }}
For convenience we report the essential information about the 
transformations we need, from Dwork's tables in [Ku]:
$$ \eqalign{
&  \l[ \matrix{
      \theta_9(\lambda) = \lambda^{-1} \hfill &\quad
      M_9(\vec a)=(a_1,a_1-a_3,a_1-a_2)  \hfill \cr
      h_9(\vec a, \lambda) = \m^{a_3-a_1-a_2} \lambda^{-a_1} \hfill &\quad
      N_9(\vec a)=\pmatrix{1 & 0 \cr 1 & -1 \cr}  \hfill \cr
      h_9(\vec a,\vec b, \lambda) =
               \m^{\mu_1+\mu_2-\mu_3} \lambda^{-\mu_1}\hidewidth\hfill \cr} \r.
\cr
&  \l[ \matrix{
      \theta_5(\lambda) = 1-\lambda \hfill & \qquad\qquad\qquad
      M_5(\vec a)=(a_1,a_2,a_1+a_2-a_3)  \hfill \cr
      h_5(\vec a, \lambda) = \m^{a_2}  \hfill & \qquad\qquad\qquad
      N_5(\vec a)=\pmatrix{0 & 1 \cr 1 & 0 \cr}  \hfill \cr
      h_5(\vec a,\vec b, \lambda) =
               \m^{-\mu_2} \hidewidth\hfill \cr}\r.
\cr
&  \l[ \matrix{
      \theta_{11}(\lambda) = (1-\lambda)^{-1} \hfill &
      M_{11}(\vec a)=(a_1,a_3-a_2,a_1-a_2)  \hfill \cr
      h_{11}(\vec a, \lambda) = \m^{a_3-a_2} (1-\lambda)^{-a_1} \hfill &
      N_{11}(\vec a)=\pmatrix{1 & -1 \cr 1 & 0 \cr}  \hfill \cr
      h_{11}(\vec a,\vec b, \lambda) =
            \m^{\mu_2-\mu_3} {{(1-\lambda^p)^{b_1}} \over {(1-\lambda)^{a_1}}}
          \hidewidth\hfill \cr}\r.
\cr}
$$
Moreover $H_m(\vec a, \lambda)=h_m(\vec a, \lambda) N_m(\vec a)$.

  From [Ku, \S 2], for  a solution matrix $C_{\vec a}(x,\lambda)$ of the
hypergeometric system at a point $x$, then
$C_{M_m(\vec a)}(x,\theta_m(\lambda)) H_m(\vec a, \lambda)^t $ is a
solution matrix of the hypergeometric equation at $\theta_m(x)$.
So
$$
     \pmatrix{ 1 & 0 \cr 0 & \theta_9(\lambda)^{-M_9(\vec a)_3} \cr}
     U_{M_9(\vec a)}^{(0)}(\theta_9(\lambda)) H_9(\vec a, \lambda)^t
     = \pmatrix{\lambda^{-a_1} & 0 \cr 0 & \lambda^{-a_2}}
     U_{M_9(\vec a)}^{(0)}(\theta_9(\lambda)) N_9^t
$$
is a solution at infinity; comparing with the solution in \copy\riferA \
we obtain the first formula.

For the second we start with
$$ \pmatrix{ 1 & 0 \cr 0 & \theta_5(\lambda)^{-M_5(\vec a)_3} \cr}
     U_{M_5(\vec a)}^{(0)}(\theta_5(\lambda)) H_5(\vec a, \lambda)^t =
     \quad = \pmatrix{1 & 0 \cr 0 & (1-\lambda)^{a_3-a_1-a_2}}
     U_{M_5(\vec a)}^{(0)}(\theta_9(\lambda)) N_5^t  \; ,$$
a solution matrix at 1,  and compare it with the given solution at that point.
For the third we begin with
$$ \eqalign{
     &\pmatrix{ \theta_{11}(\lambda)^{-M_{11}(\vec a)_1} & 0 \cr
     0 & \theta_{11}(\lambda)^{M_{11}(\vec a)_2-M_{11}(\vec a)_3} \cr}
     U_{M_{11}(\vec a)}^{(\infty)}(\theta_{11}(\lambda))
     H_{11}(\vec a, \lambda)^t = \cr
     & \quad = \pmatrix{1 & 0 \cr 0 & (1-\lambda)^{a_3-a_1-a_2}}
     U_{M_{11}(\vec a)}^{(\infty)}(\theta_{11}(\lambda)) N_{11}^t \cr
}$$
and use the same argument.
\hfill$\square$

We now come to  the proof of our theorem.
For the second formula, we use the standard lifting of Frobenius, adapted to
the origin and to infinity, $\lambda \mapsto \lambda^p$, to  write
$$
    U_{\vec b}^{(\infty)\,\sigma} \l(\lambda^p \r)
    \gamma (\vec a, \vec b ; \lambda)^t =
    \pmatrix{ \xi_1^{(\infty)} (\vec a, \vec b ) \lambda^{\mu_1}& 0 \cr
              0 & \xi_2^{(\infty)} (\vec a, \vec b )
                 \lambda^{\mu_2} \cr}
    U_{\vec a}^{(\infty)}(\lambda )\ .
$$
Using the lemma and the following formula of [Ku \S 4]:
$$ \gamma (\vec a, \vec b ; \lambda)^t =
     h_m (\vec a, \vec b ; \lambda) N_m^\ast
     \gamma \l( M_m(\vec a), M_m(\vec b) ; \theta_m(\lambda) \r)^t N_m^t  \; ,$$
for $m=9$ and
$h_9 (\vec a, \vec b ; \lambda)=\m^{\mu_1 + \mu_2 - \mu_3} \lambda^{\mu_1}$,
we obtain
$$ \eqalign{
    & U_{M_9(\vec b)}^{(0)\,\sigma} \l( {1\over {\lambda^p}} \r)
    \gamma \l(M_9(\vec a), M_9(\vec b) ; {1 \over \lambda} \r)^t =\cr
     &\quad =\m^{\mu_1 + \mu_2 - \mu_3}
    \pmatrix{ \xi_1^{(\infty)} (\vec a, \vec b ) & 0 \cr
              0 & \xi_2^{(\infty)} (\vec a, \vec b )
                 \lambda^{\mu_2 - \mu_1} \cr}
    U_{M_9(\vec a)}^{(0)} \l({1 \over \lambda} \r) \ .\cr}
\leqno{(\infty^\prime)}
$$
On the other side, from the solution at the origin, by using $M_9(\vec a)$
instead of $\vec a$ and $1/\lambda$ instead of $\lambda$, we can write:
$$ \eqalign{
    & U_{M_9(\vec b)}^{(0)\,\sigma} \l( {1\over {\lambda^p}} \r)
    \gamma \l(M_9(\vec a), M_9(\vec b) ; {1 \over \lambda} \r)^t = \cr
    &\quad = \pmatrix{ \xi_1^{(0)} (M_9(\vec a), M_9(\vec b) ) & 0 \cr
              0 & \xi_2^{(0)} (M_9(\vec a), M_9(\vec b) )
                 \lambda^{-M_9(\vec \mu)_3} \cr}
    U_{M_9(\vec a)}^{(0)} \l({1 \over \lambda} \r) \ .\cr}
\leqno{(\infty^{\prime\prime})}
$$
Comparing ($\infty^\prime$) with ($\infty^{\prime\prime}$),
and defining $M^{(\infty)} = M_9$, we conclude.

For the first formula of the theorem, we consider the lifting of Frobenius,
adapted to one and infinity,
$\lambda \mapsto \phi (\lambda) = 1-(1-\lambda)^p$,
and we write
$$
    U_{\vec b}^{(1)\,\sigma} \l( \phi(\lambda ) \r)
    \gamma^{(\phi)} (\vec a, \vec b ; \lambda)^t =
    \pmatrix{ \xi_1^{(1)} (\vec a, \vec b ) & 0 \cr
              0 & \xi_2^{(1)} (\vec a, \vec b )
                (1-\lambda)^{\mu_1 + \mu_2 - \mu_3} \cr}
    U_{\vec a}^{(1)}(\lambda )\ .
$$
By the lemma, and  replacing $(1-\lambda)^{-1}$ by
$\lambda$, we have
$$
    U_{M_{11}(\vec b)}^{(\infty)\,\sigma} \l( \lambda^p \r) N_{11}^t
    \gamma^{(\phi)} \l(\vec a, \vec b ; 1-{1\over\lambda} \r)^t =
    \pmatrix{ \xi_1^{(1)} (\vec a, \vec b ) & 0 \cr
              0 & \xi_2^{(1)} (\vec a, \vec b )
                \lambda^{\mu_3 - \mu_1 - \mu_2} \cr}
    U_{M_{11}(\vec a)}^{(\infty)}(\lambda ) N_{11}^t\ .
$$
Now, the transformation $\theta_7$ is  inverse to $\theta_{11}$.
We use again the formula of [Ku \S 4], specialized to
$$ \gamma^{(\phi)} (\vec a, \vec b ; \theta_7(\lambda))^t =
     h_7 (M_{11}(\vec a), M_{11}(\vec b) ; \lambda)^{-1} N_7^t
     \gamma^{(\phi)} (M_{11}(\vec a), M_{11}(\vec b) ; \lambda )^t N_7^\ast $$
with
$h_7 (M_{11}(\vec a), M_{11}(\vec b) ; \lambda)=
\m^{-M_{11}(\vec \mu)_2}\lambda^{M_{11}(\vec \mu)_1}=
\m^{\mu_3-\mu_2}\lambda^{\mu_1}$,
  to obtain
$$ \eqalign{
    & U_{M_{11}(\vec b)}^{(\infty)\,\sigma} \l( \lambda^p \r)
    \gamma^{(\phi)} (M_{11}(\vec a), M_{11}(\vec b) ; \lambda)^t =\cr
    & \quad =\m^{\mu_2 - \mu_3 }
    \pmatrix{ \xi_1^{(1)} (\vec a, \vec b )\lambda^{\mu_1} & 0 \cr
              0 & \xi_2^{(1)} (\vec a, \vec b )
                \lambda^{\mu_3 - \mu_2} \cr}
    U_{M_{11}(\vec a)}^{(\infty)}(\lambda ) \; .\cr }
    \leqno{(1^{\prime})}
$$
On the other side, using the solution at infinity,
and $M_{11}(\vec a)$
instead of $\vec a$ , we can write:
$$ \eqalign{
    & \qquad\qquad\qquad U_{M_{11}(\vec b)}^{(\infty)\,\sigma} \l( \lambda^p \r)
    \gamma^{(\phi)} (M_{11}(\vec a), M_{11}(\vec b) ; \lambda)^t =\cr
    &
    =\pmatrix{ \xi_1^{(\infty)} (M_{11}(\vec a), M_{11}(\vec b))
                     \lambda^{M_{11}(\vec\mu)_1} & 0 \cr
              0 & \xi_2^{(\infty)} (M_{11}(\vec a), M_{11}(\vec b) )
                \lambda^{M_{11}(\vec\mu)_2} \cr}
    U_{M_{11}(\vec a)}^{(\infty)}(\lambda )\; . \cr }
    \leqno{(1^{\prime\prime})}
$$
Again, comparing ($1^{\prime}$) and ($1^{\prime\prime}$) gives
the equalities, for $j=1,2$,
$$
   \xi_j^{(1)} (\vec a, \vec b ) =
     \m^{\mu_2 - \mu_3} \xi_j^{(\infty)} \l( M_{11}(\vec a), 
M_{11}(\vec b) \r) \; .
$$
Finally, using the formulae for $\xi_j^{(\infty)}$, we obtain
the results, where we define  $M^{(1)}:=M_5=M_9 \circ M_{11}$.
\hfill$\square$

\bigskip
{\bf \NN. {The unit root subcrystal.}}

We know that, under the conditions $\vec \mu \in \{ 0, 1, \dots , p-1 \}^3$
and
$$ \mu_3 < \min ( \mu_1 , \mu_2 )  \leqno{(1)}$$
or
$$ \mu_3 > \max ( \mu_1 , \mu_2 )\, ,  \leqno{(2)}$$
for any lifting Frobenius $\phi$ the matrix
$\gamma^{(\phi)} (\vec a, \vec b ; \lambda)$ is of the form
$$ \pmatrix{A & B \cr pC & pD \cr}  \leqno{(1)}$$
or
$$ \pmatrix{pA & pB \cr C & D \cr}\, ,  \leqno{(2)}$$
respectively,
where $A,B,C,D$ are analytic functions, bounded by $1$ on a domain
of the form $D(-{\vec \mu},p^{-1}) \times {\c S}_{\epsilon}$ where
$$ 
{\c S}_{\epsilon} = D(0, \epsilon^{-1}) \setminus
           \l( D(0,\epsilon) \cup D(1,\epsilon) \r)
$$
  for  $\epsilon \in (0,1)$.

Moreover, for each $\u \mu$ as before, there is a polynomial $F_{\vec 
\mu}(\lambda ) \in  \b Z [\lambda]$, of degree $p-1$,  
such that the 
region $\l| F_{\vec \mu}(\lambda ) \r| < 1$, consists of $p-1$ 
residue classes, called {\it supersingular}
with the property that, for any  hypergeometric system with 
parameters ${\u a} \in D(- {\u \mu}, p^{-1})$, on all classes 
$D(z,1^-)$
not singular nor supersingular, one has
$|A(z)|=1$ (resp. $|D(z)|=1$) in  case (1) (resp. (2)).

In  case (2) we define
$$\eqalign{
     \c T_2 &= \l\{ \vec a \in \b Z_p^3 \; | \;
               \mu_{a_3^{(i)}} >
               \max ( \mu_{a_1^{(i)}} , \mu_{a_2^{(i)}} ) \forall i \r\} \cr
     H_2 &= \l\{ \lambda \in \b P_{\b C_p}^{1\, rig} \; | \;
               \exists \vec a \in \c T_2 \ s.t. \
                \l| F_{\vec \mu}(\lambda ) \r| < 1 \r\} \cr
     \c S_2 &= \b P_{\b C_p}^{1\, rig} \setminus
               \l( {\rm singular\ classes} \cup H_2 \r) \; .\cr}
$$
So $\c S_2$ is the complement of a finite number of residue classes, 
and, for any lifting of Frobenius  $\phi: \; \c S_2 \to \c S_2$ we 
consider the map
$$
    \matrix{ \overline\phi \, :
             & \c T_2  \times \c S_2 & \longrightarrow & \c T_2  \times
\c S_2 \cr
             &(\vec a,\lambda ) & \longmapsto &(\vec a',\phi (\lambda) 
) \; .\cr}
$$
The search for the unit root $F$-subcrystal entails to write the 
bounded solution
inside  a residue class $D(z,1^-) \subset \c S_2$ as $(\eta u, u)$. 
Such a solution is an
eigenvector of Frobenius with a unit eigenvalue:
$$
    \overline\phi^{\ast} (\eta u, u) \gamma^{(\phi)}=\xi (\eta u, u)
$$
with $|\xi|=1$. So we have
$$
    \eta={{pA\overline\phi^{\ast} (\eta) + C } \over
          {pB\overline\phi^{\ast} (\eta) + D }}
$$
and this proves the analyticity of $\eta$ in $\c T_2  \times \c S_2$,
because the function
$$
\omega \longmapsto {{pA\overline\phi^{\ast} (\omega ) + C } \over
                     {pB\overline\phi^{\ast} (\omega ) + D }} \leqno{(\c M)}
$$
is a contraction of the Banach space of analytic functions bounded by $1$
on $\c T_2  \times \c S_2$.

The unit root $F$-subcrystal is then defined over $\c T_2  \times \c S_2$ by
$$
{u' \over u} = {{a_3 - a_1}\over {1-\lambda}} \eta
                   + {{a_1 + a_2 - a_3}\over {1-\lambda}} \ .
$$
Our main result is that, if $\phi$ is adapted at the singular point
$z \in \{ 0, 1, \infty \}$, then the map $\c M$ is a contraction of the
space of analytic functions bounded by $1$ on $\c T_2 \times D(z,1^-)$.
Therefore $\eta$ admits an extension inside the three singular 
classes (unless they are at the same time supersingular!)
and therefore  the unit root $F$-subcrystal also extends as a
logarithmic $F$-subcrystal of the hypergeometric system everywhere 
except on the supersingular locus.
Moreover, this logarithmic $F$-subcrystal
is not singular ({\it i.e.} it is an $F$-crystal in the usual sense)
in the class $D(z,1^-)$ for $z \in \{ 0, 1, \infty \}$
if and only if the bounded solution
is holomorphic, that is, if and only if $\l| \xi_1^{(z)} \r| =1$.

\bigskip
{\bf \NN. {The Koblitz-Diamond formula.}}

Under the hypothesis $\vec a \in \c T_2$ we have
$|\xi_1^{(0)} (\vec a, \vec a' ) |=1$, so the first row of $U_{\vec a}^{(0)}$
is bounded  by $1$, {\it i.e.} the unit $F$-subcrystal is not singular in
the class $D(0,1^-)$.
The hypothesis of the Koblitz-Diamond theorem
implies in fact that the same is
true in the class $D(1,1^-)$.
Then under these hypotheses the unit root $F$-subcrystal is  a 
crystal in the usual
sense in a region containing both classes
$D(0,1^-)$ and $D(1,1^-)$.
The bounded solution in these two classes is
$ \l( u_1^{(z)}, u_2^{(z)} \r)$ and $\eta=\ds{ u_1^{(z)} \over u_2^{(z)} }$.

Let $\phi$ be a lifting of Frobenius to the formal affine line $\hat {\b A}$,
adapted at $0$ and $1$: for example one could take
$\phi \l( {\lambda \over {1-\lambda}} \r) =
         \l( {\lambda \over {1-\lambda}} \r)^p
$.
We note that, if  $\theta$ indicates the transformation
$\lambda \mapsto \ds{\lambda \over {1-\lambda}}$, one obtains
$$ \lambda = \theta^{-1} \l( \theta (\lambda ) \r)
             = {t \over {1+t}} \circ {\lambda \over {1-\lambda}}
$$
therefore
$$
\phi(\lambda)= \theta^{-1} \l( \theta (\lambda )^p \r)
={{\lambda^p /(1-\lambda)^p}\over {1+\lambda^p /(1-\lambda)^p }}
={{\lambda^p}\over {(1-\lambda)^p +\lambda^p }}
=\lambda^p \big( 1-p\lambda P(\lambda ) \big)^{-1}
$$
where
$P(\lambda )= \sum\limits_{i=1}^{p-1}
      \ds{1\over p} \ds{p \choose i}(-\lambda)^{i-1}$,
hence $\phi(\lambda)  \in \c O \l( \hat {\b A} \r)$.
In similar way we obtain
$$\phi(1-\lambda)=(1-\lambda)^p \big( 1-p(1-\lambda) P(1-\lambda )
\big)^{-1}.$$

  From the formulas
$$ \overline\phi^{\ast} \l( u_1^{(z)}, u_2^{(z)} \r) \gamma^{(\phi)} =
     \xi_1^{(z)}  \l( u_1^{(z)}, u_2^{(z)} \r)
$$
for $z \in \{ 0,1 \}$, we obtain
$$ \overline\phi^{\ast} \l( u_1^{(z)} \r) \gamma_{12}^{(\phi)} +
     \overline\phi^{\ast} \l( u_2^{(z)} \r) \gamma_{22}^{(\phi)} =
     \xi_1^{(z)}  u_2^{(z)}
$$
and (dividing by $\overline\phi^{\ast} \l( u_2^{(z)} \r)$)
$$ \overline\phi^{\ast} \l( \eta \r) \gamma_{12}^{(\phi)} +
     \gamma_{22}^{(\phi)} =
     \xi_1^{(z)}  {{ u_2^{(z)} } \over {\overline\phi^{\ast} \l( u_2^{(z)} \r)}}
$$
for $z \in \{ 0,1 \}$.
Then, if $\phi$ is adapted to $0$ and $1$, we have the equality
$$ \xi_1^{(0)} (\vec a, \vec a' ) {a_3 \over a_3'}
     \c F^{(\phi)} (\vec a; \lambda )=
     \xi_1^{(1)} (\vec a, \vec a' ) {{a_1-a_3} \over {a_1'-a_3'}}
     \c F^{(\phi)} (a_1,a_2,a_1 + a_2 - a_3 +1; 1-\lambda )
$$
as analytic functions on
$\c T_2 \times \l( \c S_2 \cup D(0,1^-) \cup D(1,1^-) \r) $.
Now, we are dealing with  the Frobenius matrix $F(\phi)$ for an 
$F$-crystal {\it non-singular
at $0$ and $1$}, so we can forget the restriction that $\phi$ be 
adapted to the singularities, that is we can 
replace $\phi$ in the 
equality with the standard Frobenius
$\lambda \mapsto \lambda^p$. Evaluation at $\lambda = 1$ then gives
$$ \xi_1^{(0)} (\vec a, \vec a' ) {a_3 \over a_3'}
     \c F (\vec a; 1)=
     \xi_1^{(1)} (\vec a, \vec a' ) {{a_1-a_3} \over {a_1'-a_3'}}
$$
{\it i.e.}
$$ \c F (\vec a; 1)=
     {a_3' \over a_3} \cdot{{a_1-a_3} \over {a_1'-a_3'}} \cdot
     {{\xi_1^{(1)} (\vec a, \vec a' )} \over
     {\xi_1^{(0)} (\vec a, \vec a' )}} \; .
$$
This equality  is essentially the  Koblitz-Diamond formula. In fact,
for $p \vec b - \vec a = \vec \mu $, we can rewrite the r.h.s. as
$$\eqalign{
     &\quad {b_3 \over a_3} \cdot{{a_1-a_3} \over {b_1-b_3}} \cdot
     {{\xi_1^{(1)} (\vec a, \vec b )} \over
     {\xi_1^{(0)} (\vec a, \vec b )}} =\cr
     &=\m^{\mu_2} {b_3 \over a_3}\cdot {{a_1-a_3} \over {b_1-b_3}}\cdot
     {{\xi_1^{(0)} (a_1,a_2,a_1+a_2-a_3, b_1,b_2,b_1+b_2-b_3 )} \over
     {\xi_1^{(0)} (\vec a, \vec b )}} \cr
     &=\m^{\mu_2} {b_3 \over a_3} \cdot{{a_1-a_3} \over {b_1-b_3}}\cdot
     {{\gp (a_2,b_2) \gp (a_1-a_3,b_1-b_3) \gp (1+a_3,1+b_3)} \over
     {\gp (a_1+a_2-a_3+1,b_1+b_2-b_3+1) \gp (a_2,b_2) \gp 
(a_3-a_2,b_3-b_2)}} \cr
     &=\m^{\mu_2}
     {{ \gp (a_1-a_3+1,b_1-b_3+1) \gp (a_3,b_3)} \over
     {\gp (a_1+a_2-a_3+1,b_1+b_2-b_3+1) \gp (a_3-a_2,b_3-b_2)}} \cr
     &=\m^{\mu_2}
     {{ \pi^{\mu_1-\mu_3+1}\Gp (a_1-a_3+1) \pi^{\mu_3}\Gp (a_3)} \over
     {\pi^{\mu_1+\mu_2-\mu_3+1}\Gp (a_1+a_2-a_3+1)
      \pi^{\mu_3-\mu_2}\Gp (a_3-a_2)}} \cr
     &=\m^{\mu_2}
     {{ \Gp (a_1-a_3+1) \Gp (a_3)} \over
     {\Gp (a_1+a_2-a_3+1) \Gp (a_3-a_2)}} \cr
     &=\m^{\mu_2} {{\m^{\mu_1-\mu_3}} \over {\m^{\mu_1+\mu_2-\mu_3}}}
     {{ \Gp (a_3-a_2-a_1) \Gp (a_3)} \over
     {\Gp (a_3-a_1) \Gp (a_3-a_2)}} \cr
     &={{ \Gp (a_3-a_2-a_1) \Gp (a_3)} \over
     {\Gp (a_3-a_1) \Gp (a_3-a_2)}} \ .\cr}
$$

\vskip 10 pt

{\bf \N {Appendix: determination of Frobenius eigenvalues {\it via}
modular relations.}}
\vskip 5 pt

{\bf \NN. {}}
Our calculations are based on the Dwork's article [Bo] \S 3.
In particular, we recall the interplay of  the base change matrices
$ B(\vec a, \vec b;\lambda)$, for $\vec a \equiv \vec b \ \mod \ \b 
Z$, of {\it loc.cit.},
with solutions
$$ C_{\vec a}(z,\lambda) B(\vec a, \vec a + \vec u;\lambda)^t =
     B(\vec a, \vec a + \vec u;z)^t C_{\vec a + \vec u}(z,\lambda) $$
and with  Frobenius matrices
$$
     B(\vec a, \vec a + \vec u;\lambda) \gamma^{(\phi)} (\vec a,\vec b;\lambda)=
     \gamma^{(\phi)} (\vec a+ \vec u,\vec b+ \vec v;\lambda)
     B(\vec b, \vec b + \vec v;\phi(\lambda))
$$
(the second formula follows from the first using the Frobenius action
\copy\riferB \
on solutions).

{\bf \NNN. {}}
For $z =  0,1,\infty $,   we use the solution matrix of \copy\riferA
$$ C_{\vec a}(z,\lambda) =
     l_z(\lambda)^{D_z(\vec a)} U^{(z)}_{\vec a}(\lambda) $$
where $l_z(\lambda) = \lambda$, $1-\lambda$, $\lambda^{-1}$, and 

$$D_z(\vec a) = \pmatrix{0&0\cr 0& -a_3} \; , \;\;\; \pmatrix{0&0\cr 
0& a_3-a_1-a_2} \; , \;\;\;   \pmatrix{a_2&\cr 0& a_3} \; ,
$$
respectively.
We notice that
$B(\vec a, \vec a + \vec u;z)^t C_{\vec a + \vec u}(z,\lambda) $
is a solution at $z$ of the system
${{dY} \over {d\lambda}} = Y G_{\vec a + \vec u}(\lambda )$,
so that we can write
$$
    l_z(\lambda)^{D_z(\vec a)} U^{(z)}_{\vec a}(\lambda)
    B(\vec a, \vec a + \vec u;\lambda)^t =
    \Delta l_z(\lambda)^{D_z(\vec a + \vec u)}
    U^{(z)}_{\vec a + \vec u}(\lambda)
$$
with $\Delta$ diagonal and, using that $l_z(\lambda)^{D_z(\vec a)}$
is a diagonal matrix, the formula gives
$$
    U^{(z)}_{\vec a + \vec u}(\lambda) =
    \pmatrix{\alpha^{(z)}_1(\vec a, \vec u) & 0 \cr
              0 & \alpha^{(z)}_2(\vec a, \vec u) }
    l_z(\lambda)^{- D_z(\vec u)} 
    U^{(z)}_{\vec a}(\lambda)
    B(\vec a, \vec a + \vec u;\lambda)^t\ .
\leqno{(\bf \NNN)}
$$
\newbox\riferFRFO \global\setbox\riferFRFO=\copy\corrente

{\bf \NN}
  From this we make explicit calculations at the singular points
using the solutions written in \copy\riferA.

\proclaim {{\bf \NNN.} {Lemma}}. {
With the previous notation we have
the following values for the terms $\alpha^{(z)}_i$:
\endgraf\noindent
at the origin
$$
    \matrix{ \alpha^{(0)}_1(\vec a, \vec e_1)=1 \hfill
               &&\alpha^{(0)}_2(\vec a, \vec e_1)=
                 \ds{{a_1}\over{a_1-a_3}} \hfill\cr
             \alpha^{(0)}_1(\vec a, \vec e_2)=\ds{{a_3-a_2-1}\over{a_2}} \hfill
               &&\alpha^{(0)}_2(\vec a, \vec e_2)=-1 \hfill\cr
             \alpha^{(0)}_1(\vec a, \vec e_3)=\ds{{a_3+1}\over{a_3-a_2}} \hfill
               &&\alpha^{(0)}_2(\vec a, \vec e_3)=
                \ds{{a_1-a_3-1}\over{a_3-1}} \hfill\cr
}
$$
at $1$
$$
    \matrix{ \alpha^{(1)}_1(\vec a, \vec e_1)=
              \ds{{a_1+a_2-a_3+1}\over{a_1-a_3}} \hfill
               &&\alpha^{(1)}_2(\vec a, \vec e_1)=
               \ds{{a_1}\over{a_1+a_2-a_3-1}} \hfill\cr
             \alpha^{(1)}_1(\vec a, \vec e_2)=
                \ds{{a_3-a_1-a_2-1}\over{a_2}} \hfill
               &&\alpha^{(1)}_2(\vec a, \vec e_2)=
                \ds{{a_2-a_3+1}\over{a_3-a_1-a_2+1}} \hfill\cr
             \alpha^{(1)}_1(\vec a, \vec e_3)=
                \ds{{a_1-a_3-1}\over{a_1+a_2-a_3}} \hfill
               &&\alpha^{(1)}_2(\vec a, \vec e_3)=
                 \ds{{a_1+a_2-a_3-2}\over{a_3-a_2}} \hfill\cr
}
$$
and at infinity
$$
    \matrix{ \alpha^{(\infty)}_1(\vec a, \vec e_1)=
              \ds{{a_1-a_2+1}\over{a_3-a_1}}\hfill 
               &&\alpha^{(\infty)}_2(\vec a, \vec e_1)=
                 \ds{{a_1}\over{a_1-a_2-1}} \hfill\cr
             \alpha^{(\infty)}_1(\vec a, \vec e_2)=
              \ds{{a_3-a_2-1}\over{a_2-a_1}} \hfill
               &&\alpha^{(\infty)}_2(\vec a, \vec e_2)=
                \ds{{a_2-a_1+2}\over{a_2}} \hfill\cr 
             \alpha^{(\infty)}_1(\vec a, \vec e_3)=
                \ds{{a_3-a_1+1}\over{a_3-a_2}} \hfill
               &&\alpha^{(\infty)}_2(\vec a, \vec e_3)=1 \ . \hfill\cr
}
$$
}\par
{\proc {Proof. }}
The computation follows in every case the strategy of the
previous section.
For example, in the calculation of $\alpha^{(0)}_i(\vec a, \vec e_1)$,
we specialize $\vec u=\vec e_1=(1,0,0)$ so
$\vec a + \vec u=(a_1+1,a_2,a_3)$.
  From [Bo] we read
$$
B(\vec a, \vec a + \vec e_1;\lambda)^t = {1\over a_1}
\pmatrix{a_1-a_3 & (a_1-a_3){\lambda\over{\lambda-1}}\cr
           a_3-a_2 & {{a_1-\lambda(a_3-a_2)}\over{1-\lambda}}
           }
$$
and we can specialize our general formula (\copy\riferFRFO) to
$$
    U^{(0)}_{\vec a + \vec e_1}(\lambda) =
    \pmatrix{\alpha^{(0)}_1(\vec a, \vec e_1) & 0 \cr
              0 & \alpha^{(0)}_2(\vec a, \vec e_1) \cr}
    \b I_2
    U^{(0)}_{\vec a}(\lambda)
    B(\vec a, \vec a + \vec e_1;\lambda)^t
$$
(the singular parts in this case disappear).
Using the solution matrix of  \copy\riferA\
evaluated for $\lambda=0$
({\it i.e.} we consider only the constant term of the expansion 
w.r.t. $\lambda$),
we obtain:
$$
\pmatrix{a_3-a_2 & a_3 \cr 1-a_3 & 0 \cr }=
\pmatrix{\alpha^{(0)}_1(\vec a, \vec e_1) & 0 \cr
              0 & \alpha^{(0)}_2(\vec a, \vec e_1) \cr}
\pmatrix{a_3-a_2 & a_3 \cr 1-a_3 & 0 \cr }
{1\over a_1}
\pmatrix{a_1-a_3 & 0\cr a_3-a_2 & a_1}
$$
from which the results for $\alpha^{(0)}_1(\vec a, \vec e_1)$
and $\alpha^{(0)}_2(\vec a, \vec e_1)$ follow:
$$
a_1\pmatrix{a_3-a_2 & a_3 \cr 1-a_3 & 0 \cr }=
\pmatrix{\alpha^{(0)}_1(\vec a, \vec e_1) & 0 \cr
              0 & \alpha^{(0)}_2(\vec a, \vec e_1) \cr}
\pmatrix{a_1 (a_3-a_2) & a_1 a_3 \cr (1-a_3)(a_1-a_3) & 0 \cr }.
$$
In other cases a more precise evaluation of the matrix
solution is necessary.
For example, in the calculation of $\alpha^{(1)}_i(\vec a, \vec e_1)$
we can specialize our general formula (\copy\riferFRFO) to
$$
    U^{(1)}_{\vec a + \vec e_1}(\lambda) =
    \pmatrix{\alpha^{(1)}_1(\vec a, \vec e_1) & 0 \cr
              0 & \alpha^{(1)}_2(\vec a, \vec e_1) \cr}
    \pmatrix{1&0\cr 0&(1-\lambda)\cr}
    U^{(1)}_{\vec a}(\lambda)
    B(\vec a, \vec a + \vec e_1;\lambda)^t \ .
$$
In this case we  use the evaluation at $\lambda=1$ and
it is convenient to write the matrix $B(\vec a, \vec a + \vec e_1;\lambda)^t$
in the form
$$
B(\vec a, \vec a + \vec e_1;\lambda)^t = {1\over a_1}
\pmatrix{a_1-a_3 & (a_1-a_3)+(a_1-a_3){1\over{\lambda-1}}\cr
           a_3-a_2 & (a_3-a_2)+(a_1+a_2-a_3){1\over{1-\lambda}}
           } \; .
$$
We need the lowest order term in the expansion w.r.t. $(1-\lambda)$ for the
solution given in \copy\riferA.
We obtain the equality
$$
\eqalign{&\qquad
a_1\pmatrix{a_1+a_2-a_3+1 & a_1-a_3+1 \cr 0 & a_1+a_2-a_3 \cr }=\cr &=
\pmatrix{\alpha^{(1)}_1(\vec a, \vec e_1) & 0 \cr
              0 & \alpha^{(1)}_2(\vec a, \vec e_1) \cr}
\pmatrix{a_1(a_1-a_3) & a_1(a_1-a_3){a_1-a_3+1\over a_1+a_2-a_3+1}\cr
0 & (a_1+a_2-a_3)(a_1+a_2-a_3+1)}
}$$
from which the result follows.
The other cases are obtained in a similar way.
\hfill$\square$

\vskip 5pt
{\bf \NNN. {}}
We have  the obvious relations
$$
    \alpha^{(z)}_i(\vec a, \vec u + \vec v )=
    \alpha^{(z)}_i(\vec a + \vec u, \vec v )
    \alpha^{(z)}_i(\vec a , \vec u )
\qquad\hbox{ and } \qquad
\alpha^{(z)}_i(\vec a, \vec 0 )=1\ .
$$
So, if we write
$$
    \alpha^{(z)}_i(\vec a , \vec u )=
    {{\phi^{(z)}_i(\vec a + \vec u )}\over{\phi^{(z)}_i(\vec a)}}
$$
  we have the following

{\proc {Lemma. }}
$$
\eqalign{
   \phi^{(0)}_1(\vec a)&=\ds{{\G(a_3+1)}\over{\G(a_2)\G(a_3-a_2)}} \hfill
    \cr
   \phi^{(1)}_1(\vec a)&=\m^{a_2}
\ds{{\G(a_1+a_2-a_3+1)}\over{\G(a_2)\G(a_1-a_3)}}\hfill
   \cr
   \phi^{(\infty)}_1(\vec a)&= \m^{a_1+a_2+a_3}
\ds{{\G(a_1-a_2+1)}\over{\G(a_3-a_2)\G(a_1-a_3)}}
\cr
} \qquad
\eqalign{
    \phi^{(0)}_2(\vec a)&=\m^{a_2} \ds{{\G(a_1)}\over{\G(a_1-a_3)\G(a_3-1)}}
   \cr
   \phi^{(1)}_2(\vec a)&=\ds{{\G(a_1)}\over{\G(a_1+a_2-a_3-1)\G(a_3-a_2)}}
   \cr
\phi^{(\infty)}_2(\vec a)&= \m^{a_2}
\ds{{\G(a_1)}\over{\G(a_1-a_2-1)\G(a_2)}} \ . 
\cr
}
$$
{\proc {Proof. }}
This follows directly from the previous lemma, using the
functional equation
$\G(x+1)=x\G(x)$ ({\it i.e.} $\G(x-1)=\G(x)/(x-1)$) for the 
(classical) Gamma function.
\hfill$\square$

\vskip 5pt
{\bf \NN}
We now deduce  the modular properties of the Frobenius eigenvalues
$\xi_i^{(z)} (\vec a, \vec b )$ in terms of the $\alpha^{(z)}_i(\vec
a , \vec u )$.
The Frobenius action of  \copy\riferB \
adapted to a singular point,
on  solutions at the same point, gives
$$
    l_z(\phi(\lambda))^{D_z(\vec b)}
    U_{\vec b}^{(z)} (\phi(\lambda))
    \gamma^{(\phi)} (\vec a, \vec b ; \lambda)^t =
    \pmatrix{ \xi_1^{(z)} (\vec a, \vec b ) & 0 \cr
              0 & \xi_2^{(z)} (\vec a, \vec b ) \cr}
    l_z(\lambda)^{D_z(\vec a)}
    U_{\vec a}^{(z)}(\lambda )
$$
and, the singular part of solutions being diagonal,
$$
    U_{\vec b}^{(z)} (\phi(\lambda))
    \gamma^{(\phi)} (\vec a, \vec b ; \lambda)^t =
    \pmatrix{ \xi_1^{(z)} (\vec a, \vec b ) & 0 \cr
              0 & \xi_2^{(z)} (\vec a, \vec b ) \cr}
    l_z(\lambda)^{D_z(\vec a)}
    l_z(\phi(\lambda))^{-D_z(\vec b)}
    U_{\vec a}^{(z)}(\lambda ) \; .
$$
We make the substitutions $\vec a$ to $\vec a+\vec u$ and
$\vec b$ to $\vec b+\vec v$:
$$
\eqalign{&
    U_{\vec b +\vec v}^{(z)} (\phi(\lambda))
    \gamma^{(\phi)} (\vec a +\vec u, \vec b +\vec v ; \lambda)^t = \cr
&\quad
    \pmatrix{ \xi_1^{(z)} (\vec a +\vec u, \vec b +\vec v ) & 0 \cr
              0 & \xi_2^{(z)} (\vec a +\vec u, \vec b +\vec v ) \cr}
    l_z(\lambda)^{D_z(\vec a +\vec u)}
    l_z(\phi(\lambda))^{-D_z(\vec b +\vec v)}
    U_{\vec a +\vec u}^{(z)}(\lambda )  \; ,\cr
}
$$
then we apply the translation formulas
on  solutions
and  Frobenius:
$$
\eqalign{&
    \pmatrix{\alpha^{(z)}_1(\vec b, \vec v) & 0 \cr
              0 & \alpha^{(z)}_2(\vec b, \vec v) }
    l_z(\phi(\lambda))^{D_z(\vec b) - D_z(\vec b + \vec v)}
    U^{(z)}_{\vec b}(\phi(\lambda))
    B(\vec b, \vec b + \vec v;\phi(\lambda))^t \cr
&\quad
    B(\vec b, \vec b + \vec v;\phi(\lambda))^\ast
    \gamma^{(\phi)} (\vec a, \vec b ; \lambda)^t
    B(\vec a, \vec a + \vec u;\lambda)^t = \cr
&
    =\pmatrix{ \xi_1^{(z)} (\vec a +\vec u, \vec b +\vec v ) & 0 \cr
              0 & \xi_2^{(z)} (\vec a +\vec u, \vec b +\vec v ) \cr}
    l_z(\lambda)^{D_z(\vec a +\vec u)}
    l_z(\phi(\lambda))^{-D_z(\vec b +\vec v)} \cr
&\quad
    \pmatrix{\alpha^{(z)}_1(\vec a, \vec u) & 0 \cr
              0 & \alpha^{(z)}_2(\vec a, \vec u) }
    l_z(\lambda)^{D_z(\vec a) - D_z(\vec a + \vec u)}
    U^{(z)}_{\vec a}(\lambda)
    B(\vec a, \vec a + \vec u;\lambda)^t \; .\cr
}
$$
We finally  make the possible simplifications and use
   Frobenius on the left
$$
\eqalign{&\qquad
    \pmatrix{\alpha^{(z)}_1(\vec b, \vec v) & 0 \cr
              0 & \alpha^{(z)}_2(\vec b, \vec v) }
    \pmatrix{ \xi_1^{(z)} (\vec a, \vec b ) & 0 \cr
              0 & \xi_2^{(z)} (\vec a, \vec b ) \cr}
    l_z(\lambda)^{D_z(\vec a)}
    U_{\vec a}^{(z)}(\lambda ) = \cr
&
    =\pmatrix{ \xi_1^{(z)} (\vec a +\vec u, \vec b +\vec v ) & 0 \cr
              0 & \xi_2^{(z)} (\vec a +\vec u, \vec b +\vec v ) \cr}
    \pmatrix{\alpha^{(z)}_1(\vec a, \vec u) & 0 \cr
              0 & \alpha^{(z)}_2(\vec a, \vec u) }
    l_z(\lambda)^{D_z(\vec a)}
    U^{(z)}_{\vec a}(\lambda) \ .
}
$$
  From this formula we deduce
$$
\eqalign{&\quad
    \pmatrix{\alpha^{(z)}_1(\vec b, \vec v) & 0 \cr
              0 & \alpha^{(z)}_2(\vec b, \vec v) }
    \pmatrix{ \xi_1^{(z)} (\vec a, \vec b ) & 0 \cr
              0 & \xi_2^{(z)} (\vec a, \vec b ) \cr} = \cr
&=
    \pmatrix{ \xi_1^{(z)} (\vec a +\vec u, \vec b +\vec v ) & 0 \cr
              0 & \xi_2^{(z)} (\vec a +\vec u, \vec b +\vec v ) \cr}
    \pmatrix{\alpha^{(z)}_1(\vec a, \vec u) & 0 \cr
              0 & \alpha^{(z)}_2(\vec a, \vec u) } \cr
}
$$
that is
$$
     \xi_i^{(z)} (\vec a +\vec u, \vec b +\vec v ) =
     {{\alpha^{(z)}_i(\vec b, \vec v)}\over{\alpha^{(z)}_i(\vec a, \vec u)}}
     \xi_i^{(z)} (\vec a, \vec b )
$$
for $i=1,2$, which express the modular properties of
$\xi_i^{(z)} (\vec a, \vec b )$.

\proclaim {{\bf \NNN.} {Lemma}}. {
Let  $\theta (\vec a,\vec b)$ be a function defined for all $(\vec 
a,\vec b) \in D(0,1)^6$, with $ p\vec b -\vec a \in \b Z^3$. We 
assume that $\theta$  is a 
$p$-adic  analytic function of $(\vec 
a,\vec b)$, for fixed ${\u \mu} = p\vec b -\vec a \in \b Z^3$, 
subject to the modular properties
$$\theta (\vec a+\vec u,\vec b+\vec v)=
    {{\alpha(\vec b,\vec v)}\over{\alpha(\vec a,\vec u)}}\theta (\vec
a,\vec b) $$
for every $u,v \in \b Z^3$ with
$$\alpha(\vec a,\vec u)={{\phi(\vec a+\vec u)}\over{\phi(\vec a)}} $$
and
$$\phi (\vec a)= {{\G(\ell_1(\vec a))}
                    \over{\G(\ell_2(\vec a)) \G(\ell_3(\vec a))}} $$
where the $\ell_i(\vec a)$ are linear functions of $a_1,a_2,a_3$.
Then, up to a multiplicative constant, $\theta$ is of the form
$$\theta (\vec a,\vec b) =
     {{\gp(\ell_2(\vec a),\ell_2(\vec b)) \gp(\ell_3(\vec a),\ell_3(\vec b))}
     \over{\gp(\ell_1(\vec a),\ell_1(\vec b))}}\ .$$
}\par
{\proc Proof.} It is sufficient to check the modular
properties for $\theta (\vec a,\vec b)$: we have
$$
\eqalign{
&\quad{\theta (\vec a+\vec u,\vec b+\vec v)\over\theta (\vec a,\vec b)}=
{{\phi(\vec a)}\over{\phi(\vec a+\vec u)}}{{\phi(\vec b+\vec
v)}\over{\phi(\vec b)}} =
\cr &=
{{\G(\ell_1(\vec a))} \over{\G(\ell_2(\vec a)) \G(\ell_3(\vec a))}}
{{\G(\ell_2(\vec a+\vec u)) \G(\ell_3(\vec a+\vec
u))}\over{\G(\ell_1(\vec a+\vec u))}}
{{\G(\ell_1(\vec b+\vec v))} \over{\G(\ell_2(\vec b+\vec v))
\G(\ell_3(\vec b+\vec v))}}
{{\G(\ell_2(\vec b)) \G(\ell_3(\vec b))}\over{\G(\ell_1(\vec b))}}
\cr &=
{\G(\ell_1(\vec a))\G(\ell_1(\vec b)+\ell_1(\vec v))
\over\G(\ell_1(\vec a)+\ell_1(\vec u))\G(\ell_1(\vec b))}
{\G(\ell_2(\vec a)+\ell_2(\vec u))\G(\ell_2(\vec b))
\over\G(\ell_2(\vec a))\G(\ell_2(\vec b)+\ell_2(\vec v))}
{\G(\ell_3(\vec a)+\ell_3(\vec u))\G(\ell_3(\vec b))
\over\G(\ell_3(\vec a))\G(\ell_3(\vec b)+\ell_3(\vec v))}
}$$
and the result follows using the modular properties of the function $\gp$.
\hfill $\square$

\proclaim {{\bf \NNN.} {Corollary}}. {
For $\vec b \in D(0,1)^3$,
$p\vec b -\vec a \in \b Z^3$, then:
$$
\matrix{
       \xi^{(0)}_1(\vec a,\vec b)=\ds{{\gp(a_2,b_2)\gp(a_3-a_2,b_3-b_2)}
                   \over{\gp(a_3+1,b_3+1)}} \hfill\cr
               \xi^{(0)}_2(\vec a,\vec b)=\m^{\mu_2}
                \ds{{\gp(a_1-a_3,b_1-b_3)\gp(a_3-1,b_3-1)}
                     \over{\gp(a_1,b_1)}} \hfill\cr
       \xi^{(1)}_1(\vec a,\vec b)=\m^{\mu_2}
              \ds{{\gp(a_2,b_2)\gp(a_1-a_3,b_1-b_3)}
                  \over{\gp(a_1+a_2-a_3+1,b_1+b_2-b_3+1)}} \hfill\cr
               \xi^{(1)}_2(\vec a,\vec b)=
                 \ds{{\gp(a_1+a_2-a_3-1,b_1+b_2-b_3-1)\gp(a_3-a_2,b_3-b_2)}
                      \over{\gp(a_1,b_1)}} \hfill\cr
       \xi^{(\infty)}_1(\vec a,\vec b)= \m^{\mu_1+\mu_2+\mu_3}
              \ds{{\gp(a_3-a_2,b_3-b_2)\gp(a_1-a_3,b_1-b_3)}
                   \over{\gp(a_1-a_2+1,b_1-b_2+1)}}\hfill\cr
               \xi^{(\infty)}_2(\vec a,\vec b)= \m^{\mu_2}
                 \ds{{\gp(a_1-a_2-1,b_1-b_2-1)\gp(a_2,b_2)} 
                      \over{\gp(a_1,b_1)}} \hfill\cr
}
$$
up to multiplicative constants.
}\par

{\proc Proof.}
Use the previous lemma and the modular properties of the functions
$\xi^{(z)}_i$ for $i=1,2$ and $z\in\{0,1,\infty\}$.
\hfill $\square$


\bigskip
\centerline {\bf References.}
\bigskip
\bigskip

\item {[Ad]} Adolphson A.:
     ``Exponential sums and generalized hypergeometric functions, I:
Cohomology spaces and Frobenius action",
to appear in {\it Geometric Aspects of Dwork's Theory}, A. Adolphson,
F. Baldassarri, P. Berthelot, N. Katz, F. Loeser Eds., De Gruyter
2004.
\smallskip
\item {[Ba]} Baldassarri F.:
     ``Special values of symmetric hypergeometric functions",
     Trans. Amer. Math. Soc. 348 (1996), 2249-2289.
  \smallskip
\item {[BC1]}    Baldassarri F., Chiarellotto B.:
    `` Formal and $p$-adic theory of differential systems
    with logarithmic singularities depending upon parameters",
    Duke Math. J. 72 (1993), 241-300.
\smallskip
\item {[BC2]} Baldassarri F., Chiarellotto B.:
     `` On Andr\'e's transfer theorem", in {\it {\rm p}-adic Methods in Number
Theory and Algebraic Geometry}, A.Adolphson, S.Sperber, M.Tretkoff Eds.,
Contemporary Mathematics, Volume 133, AMS Publications 1992, 25-37.
\smallskip
\item {[Bo]} Dwork B.:
     ``On the Boyasky principle",
     Amer. J. Math. {\bf 105} (1983), 115-156.
\smallskip
\item{[Ch]} Christol G.: ``Un th\'eor\`eme de
transfert pour les disques singuliers
r\'eguliers", Ast\'erisque 119-120 (1984),151-168.
\smallskip
\item {[D]} Diamond J.:
     ``Hypergeometric series with a $p$-adic variable",
     Pacific J. Math. Soc. {\bf 94} (1981), 265-276.
\smallskip
\item{[DGS]} Dwork B., Gerotto G., Sullivan F.:``{\it An Introduction to
{\rm G}-functions}", Annals of Mathematical Studies 133, Princeton University
Press, Princeton N.J., 1994.
\smallskip
\item {[GHF]} Dwork B.:
     ``{\it Generalized Hypergeometric Functions.}",
     Oxford Mathematical Monographs,
     Clarendon Press, Oxford (1990), 188 pages.
\smallskip
    \item {[Ka]} Katz N.:
     ``Travaux de Dwork",
     S\'eminaire Bourbaki 24$^{\hbox{\`eme}}$ ann\'ee,
     1971/72, n$^o$ 409, Springer Lecture Notes 317, 34 pages.
\smallskip
    \item {[Ke]} Kedlaya K.S.:
     ``A $p$-adic local monodromy theorem",  to appear in Annals of
Math.
\smallskip
    \item {[KK]} Kato K.:
``Logarithmic structures of Fontaine-Illusie'' ,
Algebraic analysis, geometry, and number theory   (Baltimore, MD, 1988),
    191--224, Johns Hopkins Univ. Press, Baltimore, MD, 1989.
\smallskip
    \item {[Ko]} Koblitz N.:
     ``The hypergeometric function with $p$-adic parameters", in
     {\it Proceedings
      of the Queen's (Ontario) Number Theory Conference, 1979"}, 319-328.
\smallskip
    \item {[Ku]} Dwork B.:
     ``On Kummer's twenty-four solutions of the hypergeometric
      differential equation",
      Trans. Amer. Math. Soc. {\bf 285} (1984),  497-521.
\smallskip
\item {[LDE]} Dwork B.:
     ``{\it Lectures on $p$-adic Differential Equations.}",
     Grundlehren der mathematischen Wissenschaften 253, Springer (1982),
     310 pages.
\smallskip
\item{[Mz]} Mazur, B.: ``Frobenius and the Hodge
filtration",  Ann. of Math., Vol. 98 (1973), 58-95.
\smallskip
\item {[$p$-DE IV]} Dwork B.:
     ``On $p$-adic differential equations IV.
     Generalized hypergeometric functions as $p$-adic
     analytic functions in one variable",
     Annales scientifiques de l'\'Ecole Normale Sup\'erieure,
     $4^e$ S\'erie, t.6, fasc.3 (1973), 295-315.
\smallskip
\item {[Po]} Poole, E.G.C.:
     ``{\it Introduction to the Theory of Linear Differential Equations.}",
     Dover Publications Inc., New York (1960), 199 pages.
\smallskip
\item {[Sh]} Shiho A.:
     `Crystalline fundamental groups I - isocrystals on log crystalline
site and log convergent site",
     J. Math. Sci. Univ. Tokyo, 7:509-656, 2000.
\smallskip
\item {[Yo]} Young P.Th.: ``Ap\'ery numbers, Jacobi sums, and special
values of generalized $p$-adic hypergeometric functions",
    J. of Number Theory, {\bf 41} (1992), 231-255.

\bigskip  \bigskip 
\+  & {\it Francesco Baldassarri}  -  {\it 
Maurizio Cailotto} &\cr
\+
  & Dipartimento di Matematica   -  Universit\`a di Padova & \cr
\+  &
Via Belzoni 7 - I-35131 Padova - Italy & \cr
\+ &  e-mail: \cr
\+ & baldassa@math.unipd.it \cr 
\+ & cailotto@math.unipd.it
\cr
\vfill \eject

\end